\documentclass[letterpaper]{scrartcl}

\usepackage[USenglish]{babel}
\usepackage[utf8]{inputenc}

\usepackage{geometry}
\usepackage[short,nocomma]{optidef}
\usepackage{booktabs}
\usepackage{makecell}
\usepackage{subcaption}
\usepackage[hidelinks]{hyperref}
\usepackage{natbib}
\usepackage{amsmath,amssymb,amsthm}
\usepackage{algorithm}
\usepackage{algpseudocode}
\usepackage{appendix}
\usepackage{rotating}
\usepackage{tikz}
\usepackage{pgfplots}
\pgfplotsset{compat=1.17}
\usepackage{cleveref}
\usepackage{fullpage}
\usepackage{multirow}
\usepackage{doi}

\usepackage{authblk}

\title{The Bus Rapid Transit Investment Problem}

\author[1]{Rowan Hoogervorst}
\author[1]{Evelien van der Hurk}
\author[2]{Philine Schiewe}
\author[3,4]{Anita Sch\"{o}bel}
\author[3]{Reena Urban}

\affil[1]{DTU Management, Technical University of Denmark, Kongens Lyngby, 2800, Denmark,\linebreak \{rowho, evdh\}@dtu.dk}
\affil[2]{Department of Mathematics and Systems Analysis, Aalto University, Espoo, 02150, Finland, {philine.schiewe@aalto.fi}}
\affil[3]{Department of Mathematics, University of Kaiserslautern-Landau (RPTU), Kaiserslautern, 67663, Germany, {\{anita.schoebel, reena.urban\}@math.rptu.de}}
\affil[4]{Fraunhofer Institute of Industrial Mathematics ITWM, Kaiserslautern, 67663, Germany}
\date{}

\newcommand{\RR}{\mathbb{R}} 
\newcommand{\NN}{\mathbb{N}} 
\newcommand{\BRTproblem}{BRT investment problem}

\newcommand{\BRTfullTwo}[3]{\textup{\textsc{BRT}(}{#1}\textup{/}{#2}\textup{/}{#3}\textup{)}}
\newcommand{\BRTfullOne}[3]{\textup{\textsc{BRT*}(}{#1}\textup{/}{#2}\textup{/}{#3}\textup{)}}
\newcommand{\BRT}[2]{\BRTfullOne{#1}{#2}{\ensuremath{\vert M \vert \geq 1}}}
\newcommand{\SOC}[2]{\BRTfullOne{#1}{#2}{\ensuremath{\vert M \vert = 1}}}
\newcommand{\BRTtwo}[2]{\BRTfullTwo{#1}{#2}{\ensuremath{\vert M \vert \geq 1}}}
\newcommand{\SOCtwo}[2]{\BRTfullTwo{#1}{#2}{\ensuremath{\vert M \vert = 1}}}
\newcommand{\ConConstr}{BRT component constraint}
\newcommand{\Linear}{\textup{\textsc{Linear}}}
\newcommand{\MinImprov}{\textup{\textsc{MinImprov}}}
\newcommand{\EvenDemand}{\textup{EVEN}}
\newcommand{\LargeStationsDemand}{\textup{HUBS}}
\newcommand{\EndStationsDemand}{\textup{TERMINI}}
\newcommand{\PreserveBackslash}[1]{\let\temp=\\#1\let\\=\temp}
\newcolumntype{C}[1]{>{\PreserveBackslash\centering}p{#1}}
\definecolor{darkgreen}{rgb}{0,0.5,0}
\definecolor{violett}{RGB}{90,0,90}
\definecolor{azure}{rgb}{0.0, 0.5, 1.0}

\providecommand{\keywords}[1]{\small\textbf{\textit{Keywords:}} #1}
\providecommand{\supplements}[1]{\small\textbf{\textit{Supplemental material:}} #1}

\theoremstyle{plain}
\newtheorem{theorem}{Theorem}
\newtheorem{lemma}[theorem]{Lemma}
\theoremstyle{definition}
\newtheorem{definition}[theorem]{Definition}
\newtheorem{example}[theorem]{Example}

\usepackage{abstract}

\begin{document}

\maketitle

\begin{abstract}
	Bus Rapid Transit (BRT) systems can provide a fast and reliable service to passengers at low investment costs compared to tram, metro and train systems. 
	Therefore, they can be of great value to attract more passengers to use public transport.
	This paper thus focuses on the \BRTproblem{}: Which segments of a single bus line should be upgraded such that the number of newly attracted passengers is maximized?	
	Motivated by the construction of a new BRT line around Copenhagen, we consider a setting in which multiple parties are responsible for the financing of different segments of the line.
	As each party has a limited willingness to invest, we solve a bi-objective problem to quantify the trade-off between the number of attracted passengers and the investment budget.
	We model different problem variants:
	First, we consider two potential passenger responses to upgrades on the line.
	Second, to prevent scattered upgrades along the line, we consider different restrictions on the number of upgraded connected components on the line.
	We propose an epsilon-constraint-based algorithm to enumerate the complete set of non-dominated points and investigate the complexity of this problem. 
	Moreover, we perform extensive numerical experiments on artificial instances and a case study based on the BRT line around Copenhagen. 
	Our results show that we can generate the full Pareto front for real-life instances and that the resulting trade-off between investment budget and attracted passengers depends both on the origin-destination demand and on the passenger response to upgrades.
	Moreover, we illustrate how the generated Pareto plots can assist decision makers in selecting from a set of geographical route alternatives in our case study.
\end{abstract}

\noindent
\\
\supplements{\url{https://doi.org/10.11583/DTU.c.6805470}}

\section{Introduction}
Increasing the modal share of public transportation is widely recognized as an important path towards reducing greenhouse gas emissions, complementary to efforts to reduce the emissions of private cars \citep{messerli2019global}.
Bus Rapid Transit (BRT) lines can contribute to this goal, as they can offer an attractive service to passengers at relatively low investment costs compared to rail-based alternatives \citep{deng2011recent}.
A BRT line generally uses dedicated lanes for a large share of its route and is therefore not sensitive to delays as a result of traffic jams caused by private vehicles.
Moreover, BRT lines often get priority at crossings.
Therefore, BRT lines are characterized by higher speed, higher frequency, and higher reliability of service in comparison to traditional buses.

This paper concerns the planning of a single BRT line. 
Specifically, it poses the question of which segments of the BRT line should be upgraded to a full BRT standard and which could remain as a traditional mixed-traffic bus segment with the objective to maximize ridership given a limited willingness to invest.
An upgrade involves investments for the establishment of separate bus lanes as well as the upgrading of intersections and traffic installations to allow for priority of the BRT line.
Thus, the number and location of upgrades have a direct impact on the quality of a passenger's journey, and thereby on the expected ridership of the BRT line.
While there is a base amount of ridership independent of upgrades, we focus on the number of passengers that can be attracted additionally because of the improvements.
Considering the required investments in new infrastructure and the corresponding use of urban space, careful planning is needed to choose the final design of the line. 
This \BRTproblem{} can be seen as a substep of the network design phase in the traditional public transport planning process described in \citet{lusby2018survey}.

Our work is specifically motivated by the development of a new BRT line in the urban area of Copenhagen (Greater Copenhagen), which will connect multiple municipalities surrounding the city of Copenhagen \citep{movia2020brt}. 
Each of these municipalities is responsible for the investments required for the upgrading of segments that are within its borders.
Because these investments come out of their general budgets, which also cover other municipal expenses, municipalities must weigh the costs of upgrades against the societal benefits provided.
The willingness of the municipalities to work together towards a social optimum is though shown through the collaboration within the transport agency Movia, which is described in more detail in \Cref{sec:Copenhagen-case}. 
Due to this conflict of goals, we aim to quantify the impact of investments in this paper through constructing the Pareto front between the number of attracted passengers and the investment budget aggregated over all municipalities.
Moreover, a separate investment budget per municipality could lead to a bus line that often blends in and out of mixed traffic, which may not make passengers experience the line as very different from a traditional bus line.
Therefore, the \BRTproblem{} also includes a constraint to limit the overall number of upgraded connected BRT components.

In this paper, we formulate the \BRTproblem{} as a bi-objective mixed-integer linear program for two potential passenger responses to upgrades on the line: a linear and a threshold relation.
While an upgraded segment leads to a proportional number of newly attracted passengers under the \emph{linear passenger response}, passengers are only attracted to the BRT line in the \emph{threshold passenger response} if a minimum level of improvement is realized along their journey. 
The latter can be interpreted as a mode choice being made by a group of homogeneous passengers, where the passengers only switch to using the BRT line when it becomes their fastest alternative.
Considering these two different passenger responses leads to two different versions of the \BRTproblem{}, allowing us to analyze the impact of the passenger response on the trade-off between attracted passengers and investment budget.

The proposed model is intended to be used as a decision support tool within the planning process for a new BRT line.
While it can be applied in a setting with a global decision maker without municipalities to find a social optimum, its main application is the case of municipalities collaborating through a transport agency.
In this setting, a transport agency suggests solutions that are good for society and regard the concerns of the municipalities based on the generated Pareto curve.
In an iterative process, the municipalities can evaluate the suggestions and adjust their available budget until a satisfactory solution is found.
An alternative approach, which integrates a constraint representing the willingness to invest based on the number of attracted passengers per municipality within a single-objective setting, has been proposed in \cite{Hoogervorst2022}.

The contributions of this paper are four-fold.
First, we propose the bi-objective \BRTproblem{} with a \ConConstr{} and multiple investing municipalities for two alternative passenger responses to upgrades.
Second, we propose an $\epsilon$-constraint-based algorithm to solve the \BRTproblem{}, which can find all non-dominated points.
Third, in a theoretical analysis of this problem, we give tractable and intractable cases of the \BRTproblem{} and identify both NP-hard and polynomially solvable cases for the single-objective subproblems solved in our $\epsilon$-constraint-based algorithm.
Fourth, we perform an extensive computational study on artificial instances and realistic instances based on the Greater Copenhagen BRT line.

The remainder of the paper is structured as follows:
In \Cref{sec: literature}, we discuss the related literature.
In \Cref{sec:problems}, we define the problem formally, introduce the two different passenger responses and give corresponding bi-objective mixed-integer linear programming formulations.
\Cref{sec:theoretical-analysis} introduces the $\epsilon$-constraint-based algorithm used to solve the problem, and we theoretically analyze its complexity and the complexity of the single-objective problems solved within the algorithm.
We present computational results for artificial instances in \Cref{sec: artificial results}, where we analyze among others the impact of the passenger response, the \ConConstr{}, the demand pattern, and the budget split among the municipalities.
In \Cref{sec: case study results}, we describe the numerical results of our case study for the Greater Copenhagen BRT line.
The paper is concluded in \Cref{sec:conclusion}.

\section{Related Literature}
\label{sec: literature}

The public transportation planning process is traditionally split-up into a number of sequential planning steps, which range all the way from the strategic to the operational level \citep{lusby2018survey,Schiewe2022}.
The \BRTproblem{} is most closely related to the network design and line planning steps in this process, in which the public transport network and the lines operated on this network are determined.
An overview on the network design problem, and the models and solution methods used to solve it, is given by \citet{Laporte2000} and \citet{Laporte2019}. 
For an overview on the line planning problem, which is generally solved after the stations and infrastructure have been fixed, we refer to \citet{Sch10b}.
Moreover, we refer to \citet{Gattermann2017} for a discussion of the generation of line pools in line planning.

While the focus in transit network design has traditionally been on designing a transit network from scratch, recent work has increasingly focused on the improvement of existing public transport networks.
Specifically relevant for our work is the stream of literature focusing on adding dedicated bus lines within an existing transport network \citep{yao2012combinatorial,khoo2014biobjective,bayrak2018optimizing,tsitsokas2021dedicated}.
These papers focus on the placement of bus lanes along segments or lines in the network such to minimize the travel time of both bus and non-bus passengers.
This requires the evaluation of the passenger mode choice and the congestion caused by the placement of dedicated bus lanes.
While the \BRTproblem{} shares this core theme of upgrading bus segments, it focuses on a different objective: the trade-off between investment budget and attracted passengers.
Moreover, it focuses on the context of a single line and considers the effect of a constraint on the connectedness of upgraded segments.

Another relevant addition to the network design problem is the consideration of multiple investing parties.
While it is typically assumed that all investment decisions are made by one central authority, \cite{Wang2017} consider local authorities that can only make upgrade decisions for their own subgraphs, i.e., parts of the network.
In a game-theoretic setting, they formulate the interaction of the local authorities among others in a cooperative, competitive and chronological way.
Here, the aim of the local authorities is to minimize the travel time by increasing the capacity of edges under a budget constraint.
In the \BRTproblem{}, we take into account the effect of multiple municipalities through separate municipality budgets and through investigating different budget splits.
Our setting differs in considering a bi-objective problem on a single line and through the addition of the \ConConstr{}.

The underlying mathematical structure of the \BRTproblem{} also shows similarities to the more general network improvement problem.
This problem consists of choosing edges (and nodes) in a network to be upgraded while minimizing costs or satisfying budget constraints \citep{Krumke1998,Zhang2004,BaldomeroNaranjo2022}.
The problem has seen applications, e.g., in the area of road network optimization, where restricted resources can be used to upgrade edges in order to minimize the travel time between certain source-destination pairs \citep{Lin2015} or where roads can be upgraded to all-weather roads to improve the accessibility of health services \citep{Murawski2009}.
The \BRTproblem{} differs from the network improvement problem through being bi-objective and through the consideration of the \ConConstr{}.
Moreover, one of our passenger responses depends in a non-linear way on the realized improvements. 

Summarizing, the \BRTproblem{} introduced here contributes to the literature by focusing on the simultaneous consideration of the number of newly attracted passengers and the investment budget, separate municipality budgets and a \ConConstr{} within the context of upgrading one bus line.
Note that a special case of the \BRTproblem{} has been introduced in the ATMOS conference paper by \cite{Hoogervorst2022}, which looked at the single-objective problem of choosing segments to upgrade under a budget limit per investing municipality.

\section{Problem and Model Formulation}
\label{sec:problems}
In this section, we give a formal definition of the \BRTproblem{}, in which one objective reflects the passenger response and the other the investment budget.
We introduce two different passenger responses, namely \Linear{} and \MinImprov{}, and show the difference between the investment budget and the investment costs.
Finally, we provide a bi-objective mixed-integer linear programming formulation and prove its correctness.

\subsection{Problem Definition}
\label{sec:basic-setting}

The \BRTproblem{} models the selection of upgrades along a bus line.
We denote the bus line by a linear graph $(V,E)$, where the nodes $V$ represent the stations along the line and the edges $E $ denote the segments between the stations.
Upgrading a segment results in a BRT segment where the vehicles of the BRT line can operate independent from other transportation modes. 
We denote the costs of upgrading a segment $e \in E$ by $c_e \in \NN_{\geq 1}$, which encompasses all costs related to creating the infrastructure for the BRT segment.

We consider a BRT line that crosses \emph{multiple municipalities}, each of them being responsible for investments in their respective parts of the line.
We denote the set of municipalities by $M$ and let $E_m \subseteq E$ denote the set of segments within municipality $m \in M$.
We assume that the sets $E_m$ contain consecutive segments and are pairwise disjoint, which can often be achieved by splitting the segments at the borders of the municipalities.
Furthermore, we suppose that each municipality is allocated a fixed budget share $b_m$ of the (total) investment budget. 

We additionally include a \emph{\ConConstr{}} that limits the number of disjoint sequences of upgraded segments. 
We denote the maximum number of disjoint sequences by $Z \in \NN_{\geq 1}$. 
As a result of the different municipalities, each having its own budget limit, the upgraded segments might become spread-out over the BRT line without such a constraint. 
Passengers may experience such a line that constantly mixes in and out of blended traffic as not much different from a general bus line.
Moreover, such mixing into blended traffic might create delays, reducing the reliability of the BRT line and thus making connected upgrades more desirable. 
In addition, it might be easier from an organizational perspective to realize upgrades along several consecutive segments than on many (short) scattered segments.

The number of additional passengers that are attracted to the BRT line depends on the chosen segment upgrades.
We refer to this as the \textit{passenger response} to upgrades and let $p(F)$ denote the number of passengers that are newly attracted when the segments in $F \subseteq E$ are upgraded.
We evaluate two possible passenger responses: a \Linear{} passenger response in which the number of attracted passengers scales relatively to the improvement achieved on the passengers' journeys and a \MinImprov{} passenger response where passengers are attracted after a certain minimum improvement is realized along their journey.
These passenger responses are defined in \Cref{sec: passenger response}.

We are now able to define the \BRTproblem{} formally:

\begin{definition}[The \BRTproblem{}]
	\label{def:problem}
	Given are
	\begin{itemize}
		\item[] \textbf{Infrastructure:}
		\item a linear graph $(V,E)$, where $V=\{1,\ldots,n\}$ for $n\in \NN_{\geq 1}$ denotes the set of stations and $E=\bigl\{e_i=\{i,i+1\}: i\in \{1,\ldots,n-1\}\bigr\}$ the set of segments between the stations,
		\item upgrade costs $c_e\in \NN_{\geq 1}$ for all $e\in E$,
		\item an upper bound $Z\in \NN_{\geq 1}$ on the number of BRT components,
		\item[] \textbf{Municipalities:}
		\item a set of municipalities $M$,
		\item a set of consecutive segments $E_m \subseteq E$ for all $m\in M$ such that $\bigcup_{m\in M} E_m = E$ and the sets $E_m$ are pairwise disjoint,
		\item a budget share $b_m \in \RR_{>0}$ for all $m\in M$ such that $\sum_{m\in M} b_m=1$,
		\item[] \textbf{Passenger Response:}
		\item a function $p\colon 2^E \to \RR_{\geq 0}$ that determines the number of newly attracted passengers, i.e., there are $p(F)$ newly attracted passengers when upgrading the segments in $F \subseteq E$.
	\end{itemize}
	
	The aim is to determine combinations $(F,v)$ of upgraded segments~$F$ and an investment budget~$v$ that
	\begin{alignat*}{2}
		& \max p(F) &\qquad&\text{(maximize the number of newly attracted passengers)} \\
		& \min v	&&\text{(minimize the investment budget)}
	\end{alignat*}
	and satisfy the following constraints:
	\begin{itemize}
		\item The \emph{budget constraints}	
		\begin{equation*}
			\sum\limits_{e\in F \cap E_m} c_e \leq b_m v \text{ for all } m\in M \label{eq:budget1}\\
		\end{equation*} 
		restrict the investment of each municipality, where $b_m v$ is the budget of municipality $m\in M$.
		\item The \emph{\ConConstr{}} restricts the subgraph $G[F]$ induced by the segment set $F$, i.e., the subgraph of $G$ containing all edges in $F$ and their incident nodes, to have at most $Z$ connected components. 
	\end{itemize}
\end{definition}
In order to simplify notation, we call the connected components of $G[F]$ the \textit{BRT components} of $F$.
Hence, the \ConConstr{} limits the number of BRT components of $F$ to at most $Z$.

In the remainder of the paper, we are interested in finding the efficient solutions $(F,v)$ that constitute the Pareto front with respect to the number of newly attracted passengers and the investment budget:

\begin{definition}[Efficient solution, non-dominated point and Pareto front]
	\label{def:bi-objective-efficient}
	Let an instance of the \BRTproblem{} be given.
	A feasible solution $F\subseteq E$, $v\in \RR_{\geq 0}$ is called \emph{efficient} and its objective value $(p(F),v)$ is called \emph{non-dominated} if there does not exist another feasible solution $F'\subseteq E$, $v'\in \RR_{\geq 0}$ with objective value $(p(F'),v')$ such that $p(F')\geq p(F)$, $v'\leq v$ and at least one inequality holding strictly.
	The set of all non-dominated points is also called the \emph{Pareto front}.
\end{definition}

\subsection{Objective Functions Reflecting the Passenger Response} \label{sec: passenger response}
It remains to define the passenger response functions.
We model the passenger demand along the line by a set of origin-destination (OD) pairs $D\subseteq \{(i,j): i,j\in V, i\neq j\}$ that start and end at the stations of the line.
As we consider a single line, each OD pair $d \in D$ corresponds to a unique travel path $W_d \subseteq E$ along the line.
We assume that the number of potential passengers $a_d \in \NN_{\geq 1}$, who would like to travel along each OD pair $d \in D$ in case the full set of segments is upgraded, is known.
Such an estimate could follow, e.g., from a traffic study in which all segments are assumed to be upgraded. 
Moreover, passengers benefit from the infrastructure improvement $u_e\in \RR_{> 0}$ resulting from upgrading segment $e \in E$.
This improvement encompasses the reduction in travel time due to upgrading the segment, but it could, e.g., also represent the improved reliability as a result of the new BRT segment.

\begin{figure}
	\begin{subfigure}[t]{0.49\textwidth}
		\centering
		\begin{tikzpicture}
			\draw[->] (-0.2,0)-- (4.4,0); 
			\node () at (2.2,-1.3) {achieved improvement}; 
			\node[below] () at (0,-0.2) {0}; 
			\draw (4,0.1) -- (4,-0.1) node[below] {$\sum\limits_{e \in W_d} u_e$}; 
			
			\draw[->] (0,-0.2)-- (0,3.4); 
			\node[rotate=90] () at (-1.3,1.7) {attracted passengers}; 
			\draw (0.1,3) -- (-0.1,3) node[left] {$a_d$}; 
			
			\draw[dotted] (4,0) -- (4,3);
			\draw[dotted] (0,3) -- (4,3);
			
			\draw[very thick, red] (0,0) -- (4,3); 
			\draw[fill=red] (0,0) circle (0.05); 
			\draw[fill=red] (4,3) circle (0.05);
		\end{tikzpicture}
		\caption{\Linear{}}
	\end{subfigure}
	\begin{subfigure}[t]{0.49\textwidth}
		\centering
		\begin{tikzpicture}
			\draw[->] (-0.2,0)-- (4.4,0); 
			\node () at (2.2,-1.3) {achieved improvement};  
			\node[below] () at (0,-0.2) {0}; 
			\draw (4,0.1) -- (4,-0.1) node[below] {$\sum\limits_{e \in W_d} u_e$}; 
			\draw (2.5,0.1) -- (2.5,-0.1) node[below] {$L_d$}; 
			
			\draw[->] (0,-0.2)-- (0,3.4); 
			\node[rotate=90] () at (-1.3,1.7) {attracted passengers};  
			\draw (0.1,3) -- (-0.1,3) node[left] {$a_d$}; 
			
			\draw[dotted] (4,0) -- (4,3);
			\draw[dotted] (0,3) -- (4,3);
			\draw[dotted] (2.5,0) -- (2.5,3);
			
			\draw[very thick, red] (0,0) -- (2.5,0); 
			\draw[fill=red] (0,0) circle (0.05); 
			\draw[fill=white] (2.5,0) circle (0.05);
			
			\draw[very thick, red] (2.5,3) -- (4,3); 
			\draw[fill=red] (2.5,3) circle (0.05); 
			\draw[fill=red] (4,3) circle (0.05);
		\end{tikzpicture}
		\caption{\MinImprov{}}
	\end{subfigure}
	
	\caption{Illustration of the passenger responses \Linear{} and \MinImprov{}.}
	\label{fig: passenger response illustration}
\end{figure}
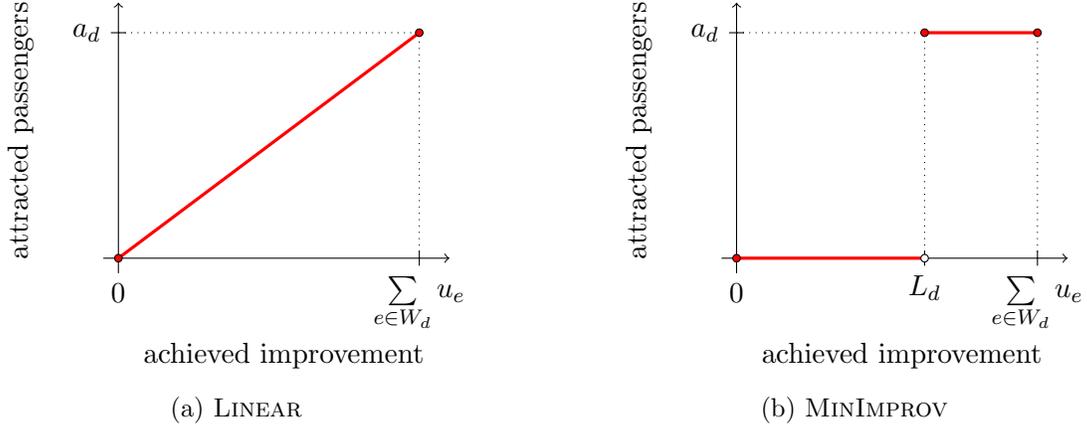

The passenger responses \Linear{} and \MinImprov{} determine the number of newly attracted passengers for each OD pair $d \in D$ based on the passenger potential $a_d$ and the infrastructure improvement realized along the path $W_d$.
These two passenger responses are illustrated in \Cref{fig: passenger response illustration}.
The \Linear{} passenger response leads to a number of attracted passengers that is proportional to the infrastructure improvement realized, i.e., realizing $x \%$ of the potential improvement leads to $x \%$ of the potential passengers being attracted.
The \MinImprov{} passenger response instead relies on a threshold $L_d \in \RR_{\geq 0}$, which represents the point at which passengers switch over to the BRT line.
An infrastructure improvement below this threshold leads to no passengers being attracted, while all potential passengers are attracted if the realized infrastructure improvement exceeds the threshold.

We now formally define the \Linear{} and \MinImprov{} passenger response:

\begin{definition}[Objective functions]
	\label{def:objective}	
	Given are
	\begin{itemize}
		\item a set of OD pairs $D\subseteq \{(i,j): i,j\in V, i\neq j\}$ with unique paths $W_d \subseteq E$ from $i$ to $j$ along the line for all $d=(i,j)\in D$,
		\item a number of potential passengers $a_d\in \NN_{\geq 1}$ for all $d\in D$,
		\item infrastructure improvements $u_e\in \RR_{>0}$ for all $e\in E$, 
	\end{itemize}
	and additionally for the \MinImprov{} passenger response:
	\begin{itemize}
		\item an improvement threshold level $L_d$ for each $d \in D$ with $L_d \leq \sum_{e\in W_d} u_e$. 
	\end{itemize}
	Let $F\subseteq E$ be the set of upgraded segments, and let an OD pair $d\in D$ be given.\\
	In \Linear{}, the number of newly attracted passengers of OD pair $d \in D$ is determined by 
	\[p_d(F) \coloneqq \frac{\sum_{e\in F\cap W_d} u_e}{\sum_{e'\in W_d}u_{e'}} \cdot a_d.\]
	In \MinImprov{}, the number of newly attracted passengers of OD pair $d \in D$ is determined by
	\[p_d(F)\coloneqq \begin{cases}
		a_d	&\text{if } L_{d} \leq \sum_{e\in F\cap W_{d}} u_e,\\
		0	&\text{otherwise}.
	\end{cases}\]
	Hence, the total number of newly attracted passengers dependent on the set of upgraded segments is given by $p\colon 2^E \to \RR_{\geq 0},\ F\mapsto \sum_{d\in D} p_d(F)$.
\end{definition}

An example of both passenger responses as well as the notation introduced in \Cref{def:problem} and \Cref{def:objective} is given in \Cref{ex:brt}.

\begin{example}	\label{ex:brt}
	Consider the example instance given in \Cref{fig:example-brt}. 
	The graph $(V,E)$ with five nodes is given at the bottom with costs $c_e$ and infrastructure improvements $u_e$ below the edges. 
	The red, dashed segments belong to municipality $m_1$ while the blue, solid segments belong to municipality $m_2$.
	The bold edges form the set $F$ of segments to be upgraded.
	Three OD pairs are given above, where the line width corresponds to the number of potential passengers $a_d$.
	
	In this example, municipality~$m_1$ invests 12 and municipality~$m_2$ invests 4.
	Because both upgraded segments in $F$ are next to each other, $F$ has only one BRT component, i.e., $F$ satisfies the \ConConstr{} for any $Z\geq 1$.
	\Cref{tab:example-brt} shows the infrastructure improvements for each OD pair as well as the number of newly attracted passengers $p_d(F)$ for \Linear{} and \MinImprov{}.
	
	\begin{figure}
		\centering
		\begin{tikzpicture}[scale=0.325, thick]
			\tikzset{
				knoten/.style={draw,circle, 
					inner sep=0, minimum size = 5pt,font=\scriptsize}
			}
			\tikzset{
				every node/.style={font=\scriptsize}
			}
			\tikzset{
				label_right/.style={font=\scriptsize, align=left, anchor=west}
			}
			\node[knoten] (v1) at (0,0) {};
			\node[knoten] (v2) at (5,0) {};
			\node[knoten] (v3) at (17,0) {};
			\node[knoten] (v4) at (21,0) {};
			\node[knoten] (v5) at (27,0) {};
			\path (v1) -- (v2) node[midway,below] {(5;4)};
			\path (v2) -- (v3) node[midway,below] {(12;11)};
			\path (v3) -- (v4) node[midway,below] {(4;4)};
			\path (v4) -- (v5) node[midway,below] {(6;5)};
			
			\draw[red, dashed, line width=3] (v2)--(v3);
			\draw[azure, line width=3] (v3)--(v4);		
			
			\node[label_right] at (28.5,-0.65) {($c_e$;$u_e$)};	
			\node[label_right] at (28.5,4) {};
			\node at (31.5,4) {$a_d$};
			\node at (34,4) {$L_d$};
			
			\draw[|-|] (0,1) --(27,1);
			\draw[line width=2] (0,1) --(27,1);
			\node[label_right] at (28.5,1) {$d_3$};	
			\node at (31.5,1) {$200$};		
			\node at (34,1) {$18$};
			
			\draw[|-|] (5,2) --(27,2);
			\draw[line width=2] (5,2) --(27,2);
			\node[label_right] at (28.5,2) {$d_2$};	
			\node at (31.5,2) {$200$};
			\node at (34,2) {$15$};
			
			\draw[|-|] (17,3) --(21,3);
			\draw[line width=1] (17,3) --(21,3);
			\node[label_right] at (28.5,3) {$d_1$};
			\node at (31.5,3) {$100$};
			\node at (34,3) {$3$};
			
			\draw[red, dashed] (v1)--(v2)--(v3);
			\draw[azure] (v3)--(v4)--(v5);		
		\end{tikzpicture}
		\caption{Example instance for the \BRTproblem{}.}\label{fig:example-brt}
	\end{figure}
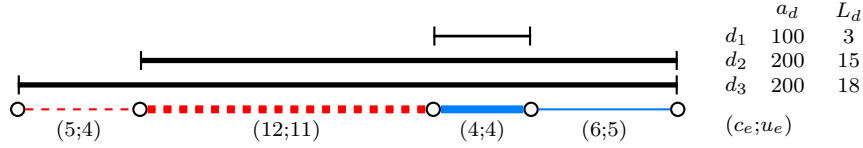
	
	\begin{table}
		\centering
		\caption{Infrastructure improvements and number of attracted passengers per OD pair for the example instance in \Cref{ex:brt}.}
		\label{tab:example-brt}
		\begin{tabular}{ccccccc}
			\toprule
			&&&&& \multicolumn{2}{c}{$p_d(F)$} \\
			\cmidrule(rl){6-7}
			OD pair & $\sum\limits_{e\in W_{d}} u_e$  & $\sum\limits_{e\in F\cap W_{d}} u_e$ & $L_d$ & $a_d$ & \Linear{} & \MinImprov{}\\
			\midrule
			$d_1$ & \phantom{2}4 & \phantom{1}4 & \phantom{1}3 & 100 & 100 & 100\\
			$d_2$ & 20 & 15 & 15 & 200 & 150 & 200\\
			$d_3$ & 24 & 15 & 18 & 200 & 125 & \phantom{10}0\\
			\bottomrule
		\end{tabular}
	\end{table}
\end{example}

\subsection{Evaluating the Investment}
\label{sec:eval_investment}

An efficient solution $(F,v)$ and its objective value $(p(F),v)$ to the \BRTproblem{} represent the set of upgraded segments, the number of newly attracted passengers and the investment budget.
For a given set of upgraded segments $F$, the \emph{investment budget} $v$ is the minimum budget such that all budget constraints are satisfied, i.e.,
\begin{equation*}
	v = \min\left\{ v' \in \RR: \sum\limits_{e\in F \cap E_m} c_e \leq b_m v' \text{ for all } m\in M\right\}.
\end{equation*}
For practical applications, however, the \emph{investment costs} $c(F)$ given as
\begin{equation*}
	c(F) \coloneqq \sum_{e\in F} c_e,
\end{equation*}
which state the actual costs incurred by upgrading the segments in $F$, are another important figure.
Because of the budget split among the municipalities based on the budget shares, for a fixed set of upgraded segments~$F$, the investment costs~$c(F)$ can be less than the available investment budget~$v$.

By solving the \BRTproblem{}, we obtain the Pareto front with respect to the investment budget.
It is not immediately clear if this Pareto front overlaps with the one where the investment costs~$c(F)$ constitute the second objective function.
We show that both Pareto fronts coincide when there is only a single municipality, i.e., $\vert M \vert = 1$, see \Cref{cor:budget-costs}.
However, this is generally not the case when there are multiple municipalities, which we illustrate with a counterexample in \Cref{ex:not-actual-investments}.

\begin{lemma}\label{cor:budget-costs}
	If $|M|=1$, then the \BRTproblem{} and the problem
	\begin{equation}
		\begin{aligned}
			\label{ip:costs}
			&\max p(F) \\
			&\min c(F) \\
			&\,\textup{s.t. there are at most } Z \textup{ BRT components},
		\end{aligned}
	\end{equation}
	where we minimize the investment costs instead of the investment budget, are equivalent in the sense that for every efficient solution of one problem there is an efficient solution of the other problem with the same objective value.
	In particular, in this case, the sets of non-dominated points coincide.
\end{lemma}
\begin{proof}
	Let $(F,v)$ be an efficient solution to the \BRTproblem{} with its corresponding non-dominated point $(p(F),v)$.
	Because $|M|=1$, the budget constraint reduces to $c(F)\leq v$.
	Because $(F,v)$ is efficient, the constraint needs to hold with equality, i.e., $c(F)=v$.
	We show that $F$ is efficient and $(p(F),v)=(p(F),c(F))$ is a non-dominated point of \eqref{ip:costs}.
	Assume that it is not efficient.
	Then there is some $F'$ such that $p(F')\geq p(F)$ and $c(F')\leq c(F)$ and at least one inequality holding strictly.
	In both cases, we have a contradiction to $(F,v)$ being efficient because the solution $(F',c(F'))$ would dominate $(F,v)$.
	
	Now let $F$ be an efficient solution to \eqref{ip:costs} with its corresponding non-dominated point $(p(F),c(F))$.
	We set $v\coloneqq c(F)$.
	Assume that $(F,v)$ is not an efficient solution to the \BRTproblem{}.
	Then there is some $(F',v')$ such that $p(F')\geq p(F)$ and $v'\leq v$ and at least one inequality holding strictly.
	Again, we have a contradiction to $F$ being efficient because the solution $F'$ would dominate~$F$ because $c(F')\leq v'\leq v=c(F)$.
\end{proof}

\begin{example} \label{ex:not-actual-investments}
	Consider the instance given in \Cref{fig: example 2 instance} with municipalities $M = \{m_1, m_2\}$ and corresponding segments $E_{m_1} = \{e_1\}$ (red, dashed) and $E_{m_2} = \{e_2\}$ (blue, solid).
	Moreover, consider a budget split in which municipality $m_1$ gets two-third and municipality $m_2$ one-third of the investment budget, i.e., $b_{m_1} = \frac{2}{3}$ and $b_{m_2} = \frac{1}{3}$, and in which there can be arbitrarily many BRT components, i.e., $Z=\infty$.
	
	For the \BRTproblem{}, the set of non-dominated points is $\{(3,3), (0,0)\}$ for the \Linear{} as well as for the \MinImprov{} passenger response due to the budget constraints.
	When considering problem~\eqref{ip:costs} (without the \ConConstr{}), it can be found that the non-dominated points are given by $\{(3,3), (2,2), (1,1), (0,0)\}$ for the \Linear{} passenger response and by $\{(3,2), (2,1), (0,0)\}$ for the \MinImprov{} passenger response.
	Comparing the sets of non-dominated points where the second objective is once the investment budget and once the investment costs, we see that they do not coincide, neither for the \Linear{} nor for the \MinImprov{} passenger response.
	One does not even need to be contained in the other.
	
	\begin{figure}[htbp]
		\centering
		\begin{tikzpicture}[scale=0.325, thick]
			\tikzset{
				knoten/.style={draw,circle, 
					inner sep=0, minimum size = 5pt,font=\scriptsize}
			}
			\tikzset{
				every node/.style={font=\scriptsize}
			}
			\tikzset{
				label_right/.style={font=\scriptsize, align=left, anchor=west}
			}
			
			\node[knoten] (v1) at (0,0) {};
			\node[knoten] (v2) at (10,0) {};
			\node[knoten] (v3) at (15,0) {};
			\path (v1) -- (v2) node[midway,below] {(2;1)};
			\path (v2) -- (v3) node[midway,below] {(1;1)};
			
			\draw[red, dashed] (v1)--(v2);
			\draw[azure] (v2)--(v3);		
			
			\node[label_right] at (17.5,-0.65) {($c_e$;$u_e$)};	
			\node[label_right] at (17.5,4) {};
			\node at (20.5,3) {$a_d$};
			\node at (23,3) {$L_d$};
			
			\draw[|-|] (0,1) --(15,1);
			\draw[line width=2] (0,1) --(15,1);
			\node[label_right] at (17.5,1) {$d_2$};	
			\node at (20.5,1) {2};		
			\node at (23,1) {$1$};
			
			\draw[|-|] (0,2) --(10,2);
			\draw[line width=1] (0,2) --(10,2);
			\node[label_right] at (17.5,2) {$d_1$};	
			\node at (20.5,2) {$1$};
			\node at (23,2) {$1$};
		\end{tikzpicture}
		
		\caption{Instance for \Cref{ex:not-actual-investments} with municipality $m_1$ containing segment $e_1$ (red, dashed) and municipality $m_2$ containing segment $e_2$ (blue, solid). }
		\label{fig: example 2 instance}
	\end{figure}
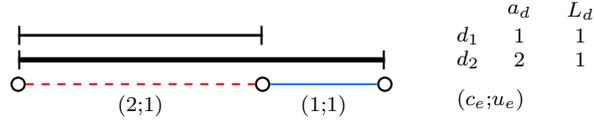	 
\end{example}

The idea of the \BRTproblem{} using the investment budget $v$ as an objective function is that the municipality budgets are relative to each other, for example, depending on sociocultural, economical or political factors.
In the numerical experiments in \Cref{sec: artificial results} and in the case study in \Cref{sec: case study results}, we compute the efficient solutions and the Pareto fronts with respect to the investment budget $v$.
Because of the practical relevance, we evaluate the results, however, also with respect to the investment costs~$c(F)$.

\subsection{Problem Variants}
\label{sec:problem-variants}

In this paper, we consider several problem variants.
We use a scheduling-like notation, where each variant of the \BRTproblem{} is classified as \BRTfullTwo{$\lambda_1$}{$\lambda_2$}{$\lambda_3$} as follows:
\begin{itemize}
	\item[] $\lambda_1$: The function chosen to represent the passenger response.
	\item[] $\lambda_2$: The upper bound on the number of BRT components of the BRT line.
	\item[] $\lambda_3$: The number of municipalities that are present.
\end{itemize}
An overview of the possible values that $\lambda_1, \lambda_2, \lambda_3$ can take is given in \Cref{tab:problem-variants}.
We remark that we use the symbolic notation ``$Z=\infty$'' to denote the setting in which the \ConConstr{} is not applied, i.e., it indicates a model without a constraint limiting the number of BRT components.
In our solution method, we also encounter the single-objective version of the problem in which we maximize the passenger response $p(F)$ given a fixed budget $v$.
These single-objective variants are classified with an asterix, i.e., as \BRTfullOne{$\lambda_1$}{$\lambda_2$}{$\lambda_3$}.

\begin{table}[htbp]
	\caption{Overview of the allowed values in the classification of the problem variants.}
	\label{tab:problem-variants}
	\centering
	\begin{tabular}{lll}
		\toprule
		Parameter & Value & Explanation \\
		\midrule
		$\lambda_1$ & \makecell[l]{\Linear{} \\ \MinImprov{} \\ $\star$} & \makecell[l]{\Linear{} passenger response \\ \MinImprov{} passenger response \\ any of the passenger responses} \\
		\midrule
		$\lambda_2$ & \makecell[l]{$Z \geq 1$ \\ $Z = k$ \\ $Z = \infty$} & \makecell[l]{any limit on the number of BRT components \\ fixed upper bound $k$ on the number of BRT components\\ no limit on the number of BRT components}\\
		\midrule
		$\lambda_3$ & \makecell[l]{$\vert M \vert \geq 1$ \\ $\vert M \vert = k$} & \makecell[l]{any number of municipalities\\ fixed number $k$ of municipalities} \\
		\bottomrule
	\end{tabular}
\end{table}

\subsection{Mixed-Integer Linear Programming Formulations}
\label{sec:IP-formulation}
We now provide a bi-objective mixed-integer linear programming formulation of the \BRTproblem{}.
This formulation uses the following variables:
\begin{itemize}
	\item a binary variable $x_e \in \{0,1\}$ for all $e\in E$ that denotes whether segment $e$ is upgraded,
	\item an auxiliary binary variable $z_{i}\in\{0,1\}$ for all $i\in \{1,\ldots,n-2\}$, which has value 1 if exactly one of the segments $e_i$ and $e_{i+1}$ is upgraded, 
	\item an auxiliary binary variable $y_d \in \{0,1\}$ for all $d\in D$ which satisfies in each optimal solution that $y_d=1$ if and only if $L_d \leq \sum_{e\in F\cap W_d} u_e$ for the set $F\subseteq E$ of upgraded segments, and
	\item a continuous variable $v \in \RR_{\geq 0}$ denoting the investment budget.
\end{itemize}

We obtain the following IP formulation, which differs with respect to the passenger response:

\noindent
\begin{tabular}{p{0.5\textwidth-2\tabcolsep}p{0.5\textwidth-2\tabcolsep}}
	\midrule
	\multicolumn{2}{C{\textwidth-2\tabcolsep}}{passenger response}\\
	\midrule
	\footnotesize{$\lambda_1 = \Linear$} & \footnotesize{$\lambda_1 = \MinImprov$}\\
	\begin{equation*}
		\begin{alignedat}{1}
			&  \max \; \sum_{e\in E} \tilde{u}_{e} x_e\\
			& \min \; v\\
			& \text{ with }\tilde{u}_e \coloneqq u_e \cdot \sum_{\substack{d\in D:\\ e\in W_d}} \frac{a_d}{\sum_{e'\in W_d}u_{e'}}
		\end{alignedat}
	\end{equation*}
	&
	\begin{equation*}
		\begin{alignedat}{1}
			& \max \; \sum\limits_{d \in D} a_d y_d \\
			& \min \; v \\
			& \text{ s.t. } L_d y_d \leq \sum_{e \in W_d} u_e x_e\quad	\text{for all } d\in D
		\end{alignedat}
	\end{equation*}
\end{tabular}
\begin{tabular}{C{\textwidth-2\tabcolsep}}
	\addlinespace
	\midrule
	budget constraints\\
	\midrule
\end{tabular}
\begin{equation*}
	\sum\limits_{e\in E_m} c_e x_e \leq b_m v\quad  \text{for all } m\in M
\end{equation*}
\begin{tabular}{C{\textwidth-2\tabcolsep}}
	\addlinespace
	\midrule
	\ConConstr{}\\
	\midrule
\end{tabular}
\begin{subequations} \label{ip:conconstr} 
	\begin{alignat}{2}
		x_{e_i}-x_{e_{i+1}} &\leq z_i &\quad& \text{for all } i\in \{1,\ldots,n-2\} \label{ip:conconstr-1}\\
		x_{e_{i+1}} - x_{e_i} & \leq z_i && \text{for all } i \in \{1,\ldots,n-2\}\label{ip:conconstr-2}\\
		x_{e_1} + \sum_{i=1}^{n-2} z_i +x_{e_{n-1}} &\leq 2Z\label{ip:conconstr-3}&&
	\end{alignat}
\end{subequations}
\begin{tabular}{C{\textwidth-2\tabcolsep}}
	\midrule
	variable domains\\
	\midrule
\end{tabular}
\begin{alignat*}{2}
	x_e &\in \{0,1\} &\quad& \text{for all } e\in E\\
	z_i &\in \{0,1\} && \text{for all } i\in \{1,\ldots,n-2\}\\
	y_d &\in \{0,1\} && \text{for all } d\in D\\
	v &\in \RR_{\geq 0}. &&\\
\end{alignat*}

In these bi-objective formulations of the \BRTproblem{}, the objectives are to maximize the number of attracted passengers and to minimize the investment budget.
The number of attracted passengers is determined either according to the \Linear{} or \MinImprov{} passenger response.
Note that the objective regarding the number of attracted passengers for $\lambda_1 = \Linear$ is reformulated as
\begin{equation*}
	\sum_{d \in D} \left( a_d\cdot \frac{\sum_{e \in W_d} u_e x_e}{\sum_{e'\in W_d}u_{e'}} \right) = \sum_{e\in E} \left(\sum_{\substack{d\in D:\\ e\in W_d}} a_d \cdot \frac{u_e}{\sum_{e'\in W_d}u_{e'}} \right)x_e = \sum_{e\in E} \tilde{u}_{e} x_e.
\end{equation*}
For $\lambda_1 = \MinImprov$ a constraint is added to ensure that the variable $y_d$ only takes value 1 in case the minimum improvement $L_d$ is realized for an OD pair $d \in D$.
The remaining constraints are the same for both passenger responses.
The budget constraints determine the available budget for each municipality based on the budget shares $b_m$. 
Moreover, the \ConConstr{}s~\eqref{ip:conconstr} ensure that the number of BRT components is no larger than $Z$.

The \ConConstr{}s~\eqref{ip:conconstr} are based on the observation that it suffices to count the number of times where an upgraded segment is succeeded by a segment that is not upgraded and vice versa.
We present the idea and its correctness formally in the following lemma:

\begin{lemma}
	\label{prop:correctness-constraints}
	Let an instance of the \BRTproblem{} be given and $F\subseteq E$.
	We reflect $F$ by setting
	$x_e \coloneqq  \begin{cases}
		1 &\text{if } e\in F,\\
		0 &\text{if } e \in E \setminus F.
	\end{cases}$\\
	Then $F$ has at most $Z$ BRT components if and only if there is a vector $z\in \{0,1\}^{n-2}$ such that the \ConConstr{}s \eqref{ip:conconstr} are satisfied.
\end{lemma}
\begin{proof}
	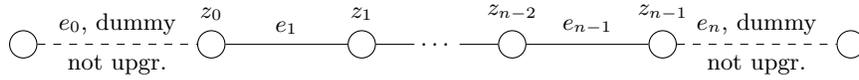
\begin{figure}[htbp]
		\centering
		\begin{tikzpicture}[font=\footnotesize]
			\node[draw, circle, fill=white] (0) at (-2.5,0) {};
			\node[draw, circle, fill=white, label=above:{$z_0$}] (1) at (0,0) {};
			\node[draw, circle, fill=white, label=above:{$z_1$}] (2) at (2,0) {};
			\node (3) at (3,0) {\dots};
			\node[draw, circle, fill=white, label=above:{$z_{n-2}$}] (4) at (4,0) {};
			\node[draw, circle, fill=white, label=above:{$z_{n-1}$}] (5) at (6,0) {};
			\node[draw, circle, fill=white] (6) at (8.5,0) {};
			\draw[dashed] (0) -- (1) node[midway, above] {$e_0$, dummy} node[midway, below] {not upgr.};
			\draw (1) -- (2) node[midway, above] {$e_1$};
			\draw (2) -- ({2.7,0}) node[midway, above] {};
			\draw ({3.3,0}) -- (4) node[midway, above] {};
			\draw (4) -- (5) node[midway, above] {$e_{n-1}$};
			\draw[dashed] (5) -- (6) node[midway, above] {$e_n$, dummy} node[midway, below] {not upgr.};
		\end{tikzpicture}
		\caption{Visualization of the \ConConstr{}s~\eqref{ip:conconstr} for \Cref{prop:correctness-constraints}.}
		\label{fig:correctness-constraints}
	\end{figure}
	Let $F \subseteq E$ be given with $K$ BRT components, i.e., $G[F]$ has $K$ connected components, denoted by $F_1,\ldots, F_K$.
	We modify the linear graph by adding dummy edges $e_0$ and $e_n$ that are not upgradable, i.e., fixed $x_{e_0}=x_{e_n}=0$, at the front and end as depicted in \Cref{fig:correctness-constraints}, and we add the binary variables $z_0\in \{0,1\}$ and $z_{n-1}\in \{0,1\}$.
	Based on that, we define $\bar{z}\in \{0,1\}^n$ by $\bar{z}_i \coloneqq \vert x_{e_i} - x_{e_{i+1}} \vert$ for all $i\in \{0,\ldots, n-1\}$.
	By definition, $\bar{z}$ is feasible for constraints~\eqref{ip:conconstr-1} and~\eqref{ip:conconstr-2}, and it has exactly $2K$ entries with value 1, namely one for each start and end of a BRT component~$F_i$, $i\in \{1,\ldots,K\}$.
	Furthermore, because $x_{e_0}=x_{e_n}=0$, we have $\bar{z}_0 = x_{e_1}$ and $\bar{z}_{n-1} = x_{e_{n-1}}$.
	
	For the first direction, let $F$ have at most $Z$ BRT components, i.e., $K\leq Z$.
	Then 
	\begin{equation*}
		x_{e_1} + \sum_{i=1}^{n-2} \bar{z}_i + x_{e_{n-1}} = \sum_{i=0}^{n-1}\bar{z}_i = 2K \leq 2Z.
	\end{equation*}
	Hence, the constraints~\eqref{ip:conconstr} are satisfied for the vector $(\bar{z}_1,\ldots, \bar{z}_{n-2})\in \{0,1\}^{n-2}$.
	
	For the second direction, we suppose that there is some $z^{\ast}\in \{0,1\}^{n-2}$ such that constraints~\eqref{ip:conconstr} hold. 
	Due to the constraints~\eqref{ip:conconstr-1} and~\eqref{ip:conconstr-2}, for all $i\in \{1,\ldots,n-2\}$, we have that $\vert x_{e_i} - x_{e_{i+1}} \vert =1$ implies $z^{\ast}_i=1$.
	Hence, $\bar{z}_i\leq z^{\ast}_i$ for all $i\in \{1,\ldots, n-2\}$.
	Then $K\leq Z$ because
	\begin{equation*}
		2K = \sum_{i=0}^{n-1}\bar{z}_i\leq x_{e_1} + \sum_{i=1}^{n-2}z^{\ast}_i +x_{e_{n-1}}\leq 2Z. \qedhere
	\end{equation*}
\end{proof}

\section{Solution Method and Theoretical Analysis}
\label{sec:theoretical-analysis}

In this section, we present an algorithm based on the $\epsilon$-constraint method to solve the \BRTproblem{} and analyze it theoretically.
In particular, we analyze the complexity of the bi-objective problem, the size of its Pareto front, the impact of the \ConConstr{} and the complexity of the single-objective problems solved within the $\epsilon$-constraint method.

\subsection{Solution Method and Tractability}
Solving the \BRTproblem{} requires computing the set of non-dominated points of an instance of \BRTtwo{$\star$}{$Z\geq 1$}.
To do so, we employ the well-known $\epsilon$-constraint method for bi-objective programming \citep{haimes1971bicriterion}, in which a series of single-objective problems is solved by placing a bound on one of the objectives.
By varying the bound $\epsilon$ over the iterations, different solutions on the Pareto front are found.

Our algorithm for solving the \BRTproblem{}, which is an adaption of the algorithm presented by \cite{Berube2009}, is given in \Cref{alg:brt-non-dominated}.
In this algorithm, we place an upper bound on the investment budget objective, meaning that we solve single-objective problems \BRT{$\star$}{$Z \geq 1$} that contain the additional constraint $v \leq B$ for varying values of $B$.
We start by finding the budget $B$ at which all segments can be upgraded, meaning that all passengers will be attracted.
In every iteration of the algorithm, we then reduce $B$ in such a way that no non-dominated points are missed.
This is repeated as long as the budget $B$ is non-negative.

While it is common to change $\epsilon$ with a fixed step size in the $\epsilon$-constraint method, such a strategy may not find all non-dominated points.
Instead, we use the integrality of the upgrade costs to identify a step size in each iteration that does not cut-off any non-dominated point.
To do so, we first identify the minimum budget at which the current solution remains feasible and the municipalities for which this minimum budget is tight.
Due to the integrality of the upgrade costs, we know that the individual budget for each such tight municipality can be reduced by~1 without cutting off any non-dominated point.
Similarly, we can reduce the budget to the next integer level for each non-tight municipality without cutting off any non-dominated point.
We then choose the step size as the minimum value that leads to a budget satisfying these conditions for each municipality.

We formally prove that our algorithm is able to find the complete set of non-dominated points in \Cref{prop:brt-non-dominated}.

\begin{algorithm}[htbp]
	\caption{Computing the non-dominated points for \BRTtwo{$\star$}{$Z\geq 1$}}\label{alg:brt-non-dominated}
	\begin{algorithmic}[1]
		\State \textbf{Input:} instance $I$ of \BRTtwo{$\star$}{$Z\geq 1$}.
		\State \textbf{Output:} set $\Gamma$ of all non-dominated points.
		\State As start values set \\$\Gamma\gets \emptyset$,\\
		$B \gets \max\limits_{m \in M} \biggl\{\frac{1}{b_m}\cdot \sum\limits_{e\in E_m} c_e\biggr\}$,\\
		$v^*\gets \max\limits_{m \in M}\biggl\{\frac{1}{b_m}\cdot \sum\limits_{e\in E_m} c_e\biggr\}$,\\
		$p^* \gets \sum_{d\in D} a_d$.
		\While{$B\geq 0$}
		\State Compute \BRT{$\star$}{$Z\geq 1$} for instance $I$ with budget $B$. 
		Let $F$ be an optimal solution and $\bar{p}$ be the optimal objective value. \label{alg:single_obj}
		\State Compute the minimum budget $\bar{v}$ such that $F$ remains feasible as \begin{equation*}
			\bar{v} \gets \max_{m \in M} \left\{\frac{1}{b_m}\cdot \sum\limits_{e\in F \cap E_m} c_e \right\}.
		\end{equation*}\label{alg:compute_v}
		\State Determine the set of municipalities $\bar{M}(\bar{v})$ for which the budget $\bar{v}$ is tight as
		\begin{equation*}
			\bar{M}(\bar{v}) \gets \left\{m\in M: \sum\limits_{e\in F \cap E_m} c_e = b_m \cdot \bar{v}\right\}. \label{alg:tight-municipalities}
		\end{equation*}
		\State Compute step width $\delta$ as 
		\begin{equation*}
			\delta \gets \min \left\{\min\limits_{m \in \bar{M}(\bar{v})}\left\{\frac{1}{b_m} \right\}, \min\limits_{m \in M \setminus \bar{M}(\bar{v})}\left\{\frac{b_m \cdot \bar{v} - \lceil b_m \cdot \bar{v} -1 \rceil}{b_m} \right\}\right\}. \label{alg:compute_delta}
		\end{equation*}
		\If{$\bar{p}<p^*$}
		\State Update $\Gamma \gets \Gamma \cup \{(p^*,v^*)\}$.
		\State Update $p^* \gets \bar{p}$.
		\EndIf
		\State Update $v^*\gets \bar{v}$.
		\State Update $B \gets \bar{v}-\delta$.
		\EndWhile
		\State Update $\Gamma \gets \Gamma \cup \{(p^*,v^*)\}$.\\
		\Return $\Gamma$
	\end{algorithmic}		
\end{algorithm}

\begin{lemma} \label{prop:brt-non-dominated}
	\Cref{alg:brt-non-dominated} computes the set of all non-dominated points for \BRTtwo{$\star$}{$Z\geq 1$}.
\end{lemma}

\begin{proof}
	\Cref{alg:brt-non-dominated} is an adaption of the algorithm presented by \cite{Berube2009}.
	We start with $p^* = \sum_{d\in D} a_d$, which is the upper limit on the number of attracted passengers that can be realized by all municipalities $m \in M$ upgrading all segments, i.e., investing $\sum_{e\in E_m} c_e$.
	This investment is possible for each municipality if the budget is set to $B=\max\limits_{m \in M}\biggl\{\frac{1}{b_m}\cdot \sum\limits_{e\in E_m} c_e\biggl\}$.
	The idea of the algorithm is to iteratively compute all non-dominated points by solving \BRT{$\star$}{$Z\geq 1$} for a budget $B$ and then reducing~$B$ by $\delta$.
	Therefore, we have to make sure that $\delta$ is small enough to not cut off a non-dominated point (see Step~1) but not arbitrarily small such that the algorithm terminates after a finite number of steps (see Step~2).
	Checking whether a lower budget yields the same number of attracted passengers eliminates weakly dominated solutions.
	
	\paragraph{Step 1: No non-dominated point is cut off.}	
	
	First note that $\bar{v}$ as computed in \cref{alg:compute_v} is actually the smallest value such that $F$ remains feasible. 
	We need to ensure that the step width $\delta$ as computed in \cref{alg:compute_delta} does not cut off solutions with $v' < \bar{v}$.
	We do so through showing that if a set of upgrades $F' \subseteq E$ is feasible for a budget $v' < \bar{v}$, it is also feasible for the budget $\bar{v} - \delta$.
	So let $F' \subseteq E$ be feasible with a corresponding minimum investment budget $v'<\bar{v}$ and let $m\in M$ be arbitrary.
	It holds that 
	\begin{equation*}
		\sum_{e\in F' \cap E_m} c_e \leq b_m \cdot v' < b_m\cdot \bar{v}.
	\end{equation*}
	Because $c_e\in \NN_{\geq 1}$ for all $e\in E$, we obtain
	\begin{equation*}
		\sum_{e\in F' \cap E_m} c_e \leq \lfloor b_m \cdot v' \rfloor \leq \lceil b_m\cdot \bar{v} -1\rceil.
	\end{equation*} 
	Now let $\delta$ be chosen as in \cref{alg:compute_delta} and update $B\coloneqq \bar{v} - \delta$.
	This means the right-hand side of the budget constraint of municipality $m$ in the next iteration is
	\begin{equation*}
		b_m \cdot (\bar{v}- \delta) \geq \begin{cases}
			b_m \cdot \left(\bar{v}-\frac{1}{b_m} \right) = b_m \cdot \bar{v} -1 = \lceil b_m \cdot \bar{v} -1 \rceil & \text{ if } m\in \bar{M}(\bar{v}),\\[5pt]
			b_m \cdot \left(\bar{v}-\frac{b_m \cdot \bar{v} - \lceil b_m \cdot \bar{v} -1 \rceil}{b_m} \right) = \lceil b_m \cdot \bar{v} -1 \rceil & \text{ if } m\in M\setminus\bar{M}(\bar{v}),
		\end{cases}
	\end{equation*} 
	by choice of $\delta$.
	Hence, the solution $F'$ with investment budget $v'< \bar{v}$ is not cut off.
	Note that this argument works for any solutions $F', F$ with corresponding investment budgets $v'$ and $\bar{v}$, respectively, with $v'<\bar{v}$.
	
	\paragraph{Step 2: The algorithm terminates.}
	To show that the algorithm terminates, we have to show that $\delta$ cannot be arbitrarily small. Remember that $\delta$ is chosen as
	\begin{equation*}
		\delta = \min \left\{\min\limits_{m \in \bar{M}(\bar{v})}\left\{\frac{1}{b_m} \right\}, \min\limits_{m \in M \setminus \bar{M}(\bar{v})}\left\{\frac{b_m \cdot \bar{v} - \lceil b_m \cdot \bar{v} -1 \rceil}{b_m} \right\}\right\}. 
	\end{equation*}	
	Thus, $\delta >0$ is guaranteed. We additionally show that there are only finitely many possible values that $\delta$ can take.
	First note that there are only $\vert M \vert$ possible values for $\frac{1}{b_m}$ independent of~$\bar{v}$. Thus, we only have to consider possible values for $\frac{b_m \cdot \bar{v} - \lceil b_m \cdot \bar{v} -1 \rceil}{b_m}$.
	By the computation of $\bar{v}$ in \cref{alg:compute_v}, we get that $\bar{M}(\bar{v}) \neq \emptyset$ and that there is an $m \in \bar{M}(\bar{v})$ such that $\bar{v} \in \{k \cdot \frac{1}{b_m}: k \in \NN_{\geq 0}\}$.
	As $\bar{v}$ is bounded from above by $\max\limits_{m \in M} \biggl\{ \frac{1}{b_m}\cdot \sum\limits_{e\in E_m} c_e \biggr\}$, we get 
	\begin{equation*}
		\bar{v} \in \bigcup_{m \in M} \left\{ k \cdot \frac{1}{b_m} \colon k \in \NN_{\geq 0} \right\} \cap \left[0, \max\limits_{m \in M} \left\{ \frac{1}{b_m}\cdot \sum\limits_{e\in E_m} c_e\right\} \right].
	\end{equation*}
	Therefore, there are only finitely many values for $\bar{v}$ in \Cref{alg:brt-non-dominated} and by extension only finitely many values for $\frac{b_m \cdot \bar{v} - \lceil b_m \cdot \bar{v} -1 \rceil}{b_m}$ and for $\delta$. 
\end{proof}

Note that the algorithm simplifies for the special case in which there is a global decision maker, i.e., for \SOCtwo{$\star$}{$Z\geq 1$}. 
In this special case, \cref{alg:compute_v,alg:compute_delta} in \cref{alg:brt-non-dominated} simplify to
\begin{align*}
	\bar{v} & \gets \sum_{e \in F} c_e,\\
	\delta & \gets 1.
\end{align*}
This means that the minimum investment budget for a given solution corresponds to the investment costs to realize it and we can always choose the step size to be equal to 1 because of the integral costs $c_e\in \NN_{\geq 1}$.
This finding relates to \Cref{cor:budget-costs}, in which we found that the Pareto front with respect to the investment costs coincides with the one for the investment budget.

To analyze the running time of \Cref{alg:brt-non-dominated}, we have to consider both the complexity of the single-objective subproblem solved in \cref{alg:single_obj} and the number of non-dominated points. 
\Cref{thm:pareto front size} shows that the Pareto front is generally intractable, meaning that it may contain an exponential number of non-dominated points. 
We give a bound on the number of non-dominated points in \Cref{lem:non-dominated-bound}.
Moreover, we show in \Cref{lem:unit_cost_2-SOC} that the number of non-dominated points is polynomial for the special case where all segment upgrade costs are equal.

\begin{theorem}\label{thm:pareto front size}
	$\BRTtwo{\star}{Z=\infty}$ is intractable, even if $\vert M \vert =1$ and $Z=\infty$. 
\end{theorem} 

\begin{proof}
	Consider an instance of $\SOCtwo{\star}{Z=\infty}$ with a graph $(V,E)$, $\vert V \vert = n$, $D\coloneqq\{d_i \coloneqq (i,i+1): i \in \{1, \ldots, n-1\}\}$, $L_{d} = 1$ for all $d\in D$ and $u_{e} =1$ for all $e\in E$.
	Set the passenger potential to $a_{d_i}\coloneqq 2^{i-1}$ and the costs to $c_{e_i} \coloneqq 2^{i-1}$ for all $i \in \{1, \ldots, n-1\}$. 
	As the paths of all OD pairs only contain one segment, upgrading a segment~$e_i\in E$ results in attracting $a_{d_i}=2^{i-1}$ passengers both for \Linear{} and \MinImprov{}.
	
	Upgrading any set $F \subseteq E$ of segments results in attracting $\sum\limits_{i \colon e_i \in F} 2^{i-1}$ passengers with investment costs and hence also investment budget of $\sum\limits_{i \colon e_i \in F} 2^{i-1}$. 
	As each number $k \in \{0, \ldots, 2^{n-1}-1\}$ can be represented by a binary representation $k = \sum\limits_{i \in F'} 2^{i}$ for some $F' \subseteq \{0, \ldots, n-1\}$, there is a solution of $\SOCtwo{\star}{Z=\infty}$ with objective value $(k,k)$ for each $k \in \{0,\ldots, 2^{n-1}-1\}$.
	Starting with the ideal point $(2^{n-1}-1, 2^{n-1}-1)$, we can easily see that all these points are non-dominated.
	Thus the set of non-dominated points has size~$2^{n-1}$, which concludes the proof.
\end{proof}

\begin{lemma}\label{lem:non-dominated-bound}
	The number of non-dominated points for $\BRTtwo{\star}{Z=\infty}$ is limited by 
	\begin{equation*}
		1 + \sum_{m \in M} \left\lfloor b_m \cdot \max_{m' \in M} \left\{\frac{1}{b_{m'}} \cdot \sum_{e \in E_{m'}} c_e \right\}\right\rfloor.
	\end{equation*}
\end{lemma}

\begin{proof}
	\Cref{alg:brt-non-dominated} computes at most one non-dominated point per iteration, i.e., per $\bar{v}$.
	From the proof of \Cref{prop:brt-non-dominated}, we know that
	\begin{equation*}
		\bar{v} \in \bigcup_{m \in M} \left\{ k \cdot \frac{1}{b_m} \colon k \in \NN_{\geq 0} \right\} \cap \left[0, \max\limits_{m \in M} \left\{ \frac{1}{b_m}\cdot \sum\limits_{e\in E_m} c_e\right\} \right].
	\end{equation*}
	Hence, the number of non-dominated points is bounded from above by
	\begin{align*}
		&\left\vert \bigcup_{m \in M} \left\{ k \cdot \frac{1}{b_m} \colon k \in \NN_{\geq 0} \right\} \cap \left[0, \max\limits_{m \in M} \left\{ \frac{1}{b_m}\cdot \sum\limits_{e\in E_m} c_e\right\} \right] \right\vert \\
		&\leq 1 + \sum_{m\in M} \left\vert \left\{ k \in \NN_{\geq 1} : k \leq b_m \cdot \max\limits_{m' \in M} \left\{ \frac{1}{b_{m'}}\cdot \sum\limits_{e\in E_{m'}} c_e \right\} \right\} \right\vert \\
		&= 1 + \sum_{m \in M} \left\lfloor b_m \cdot \max_{m' \in M} \left\{\frac{1}{b_{m'}} \cdot \sum_{e \in E_{m'}} c_e \right\}\right\rfloor. \qedhere
	\end{align*}
\end{proof}

\begin{lemma}\label{lem:unit_cost_2-SOC}
	For $\SOCtwo{\star}{Z\geq 1}$ with $c_e=c$ for all $e \in E$, there are at most $n$ non-dominated points.
\end{lemma}

\begin{proof}
	As each solution $F\subseteq E$ has investment costs $\sum_{e \in F} c_e = \vert F \vert \cdot c$, there are at most $\vert E \vert + 1= n$ different values of investment costs for feasible solutions. 
	As the sets of non-dominated points with respect to the investment costs and with respect to the investment budget coincide by \Cref{cor:budget-costs}, there are at most~$n$ non-dominated points.
\end{proof}

Next, we identify two cases for \Linear{} with only one municipality in which the Pareto front can be computed in polynomial time, see \Cref{prop:linear-connected-pass,prop:linear-connected-cost}.
The setting from \Cref{prop:linear-connected-pass} occurs for example if we consider unit infrastructure improvements $u_e=1$ for all $e\in E$ and a passenger potential that is distributed evenly over all OD pairs, i.e., $a_d=a_{d'}$ for all $d,d'\in D$.
A cost pattern as in \Cref{prop:linear-connected-cost} occurs for example if the costs are less expensive in the middle of the line but are increasingly expensive towards its ends.
Recall that $p(F) = \sum_{e\in F} \tilde{u}_e$ for all $F\subseteq E$ with $\tilde{u}_e \coloneqq u_e \cdot \sum\limits_{\substack{d\in D:\\ e\in W_d}} \frac{a_d}{\sum_{e'\in W_d}u_{e'}}$ for all $e\in E$.

\begin{lemma}\label{prop:linear-connected-pass}
	Let an instance of \SOCtwo{\Linear}{$Z=\infty$} with unit costs ${c_e\coloneqq 1}$ for all $e\in E$ be given.
	Let $e_{(i)}$, $i\in \{1,\ldots,n-1\}$ denote a sorting of the segments such that $\tilde{u}_{e_{(1)}} \geq \ldots \geq \tilde{u}_{e_{(n-1)}}$.
	\begin{enumerate}
		\item If $v \in \{0,\dots,n-1\}$ and $F= \{e_{(i)}: i\in \{1,\ldots,v\}\}$, then $(F,v)$ is an efficient solution.
		\item If there is some $\bar{i}\in\{1,\ldots,n-1\}$ such that $\tilde{u}_{e_j}\leq \tilde{u}_{e_{j'}}$ for all $j\leq j'\leq \bar{i}$ and $\tilde{u}_{e_j}\geq \tilde{u}_{e_{j'}}$ for all $\bar{i} \leq j\leq j'$, then there is some efficient solution $(F,v)$ for each non-dominated point $(p(F),v)$ such that all segments in $F$ are connected.
		\item The instance can be solved in polynomial time.
	\end{enumerate}
\end{lemma}
\begin{proof}
	\begin{enumerate}
		\item\label{prop:linear-connected-pass-1}		
		First, $(F,v)$ is feasible because $\sum_{e\in F} c_e = v$, hence, the budget constraint is satisfied.
		Second, suppose it is not efficient.
		Then it is dominated by some solution $(F',v')$.
		Assume $v'<v$, then $\vert F' \vert = c(F') \leq v' < v = c(F) = \vert F \vert$ and hence 
		\begin{equation*}
			p(F')=\sum_{e\in F'} \tilde{u}_e < \sum_{i=1}^v \tilde{u}_{e_{(i)}} = \sum_{e\in F} \tilde{u}_e = p(F).
		\end{equation*}
		Now assume $p(F')>p(F)$.
		In this case $\vert F' \vert > \vert F \vert$ because $F$ contains the $v$ segments with the highest value $\tilde{u}_e$.
		This implies $v'\geq c(F')>c(F) =v$.
		Therefore, there cannot be a solution that dominates $(F,v)$, and $(F,v)$ is efficient.
		\item\label{prop:linear-connected-pass-2} Let $(p(F'),v)$ be a non-dominated point, i.e., $v\in \{0,\ldots,n-1\}$ because the budget constraint is satisfied with equality.
		By assumption, we can suppose that $e_{(1)}=e_{\bar{i}}$.
		Because $\tilde{u}_e$ increases monotonically until $e_{\bar{i}}$ and decreases monotonically afterwards, we can assume that $e_{(2)}\in \{e_{\bar{i}-1}, e_{\bar{i}+1}\}$.
		Iteratively, we get that if $\{e_{(1)},\ldots,e_{(k)}\} = \{e_j: j\in \{l,\ldots,l+k\}\}$ for some $l\in \{1,\ldots,n-1\}$, then $e_{(k+1)}\in \{e_{l-1},e_{l+k+1}\}$.
		Hence, $(F,v)$ with $F \coloneqq \{e_{(i)}: i\in \{1,\ldots,v\}\}$ is an efficient solution because of \cref{prop:linear-connected-pass-1} and thus $p(F)=p(F')$.
		Moreover, the set $F$ is connected.
		\item\label{prop:linear-connected-pass-3} From \Cref{lem:unit_cost_2-SOC}, we know that there are at most $n$ non-dominated points, one for each investment budget value $v\in \{0,\ldots,n-1\}$.
		The sorting of the segments with respect to the values $\tilde{u}_e$ can be done in $\mathcal{O}(n^2)$. 
		Because we can find the optimal set of upgraded segments~$F$ corresponding to a fixed investment budget value $v$ as shown in \cref{prop:linear-connected-pass-1}, we can find all non-dominated points in polynomial time by iterating over the investment budget values $v\in \{0,\ldots,n-1\}$. \qedhere
	\end{enumerate}
\end{proof}

\begin{lemma}\label{prop:linear-connected-cost}
	Let an instance of \SOCtwo{\Linear}{$Z=\infty$} with ${\tilde{u}_e\coloneqq 1}$ for all $e\in E$ be given.
	Let~$e_{(i)}$, $i\in \{1,\ldots,n-1\}$ denote a sorting of the segments such that $c_{e_{(1)}} \leq \ldots \leq c_{e_{(n-1)}}$.
	\begin{enumerate}
		\item Let $k\in \{0,\ldots,n-1\}$. If $v = \sum_{i=1}^k c_{e_{(i)}}$ and $F= \{e_{(i)}: i\in \{1,\ldots,k\}\}$, then $(F,v)$ is an efficient solution.
		\item If there is some $\bar{i}\in\{1,\ldots,n-1\}$ such that $c_{e_j}\geq c_{e_{j'}}$ for all $j\leq j'\leq \bar{i}$ and $c_{e_j}\leq c_{e_{j'}}$ for all $\bar{i} \leq j\leq j'$, then there is some efficient solution $(F,v)$ for each non-dominated point $(p(F),v)$ such that all segments in $F$ are connected.
		\item The instance can be solved in polynomial time.
	\end{enumerate}
\end{lemma}
\begin{proof}
	\begin{enumerate}
		\item\label{prop:linear-connected-cost-1} First, $(F,v)$ is feasible by construction.
		Second, suppose it is not efficient.
		Then it is dominated by some solution $(F',v')$.
		Assume $p(F')>p(F)$.
		This implies $\vert F' \vert > \vert F \vert = k$, and hence 
		\begin{equation*}
			v'\geq c(F')=\sum_{e\in F'} c_e > \sum_{i=1}^k c_{e_{(i)}} = v.
		\end{equation*}
		Now assume $v'<v$.
		In this case $\vert F' \vert < \vert F \vert$ because $F$ contains the $k$ segments with the lowest costs~$c_e$ and $c(F)=v$.
		This implies $p(F')<p(F)$.
		Therefore, there cannot be a solution that dominates $(F,v)$, and $(F,v)$ is efficient.
		\item The proof is analogous to the proof for \Cref{prop:linear-connected-pass}, \cref{prop:linear-connected-pass-2}.
		\item Because there are $n$ different values for the number of newly attracted passengers, namely $p(F)=\vert F\vert \in \{0,\ldots, n-1\}$ for all $F\subseteq E$, and $c_e>0$ for all $e\in E$, there are $n$ non-dominated points, one per number of segments in $F$.
		We can sort the segments according to their costs in $\mathcal{O}(n^2)$.
		Using the formula in \cref{prop:linear-connected-cost-1}, we can find all non-dominated points in polynomial time by iterating over the number of upgraded segments $k\in \{0,\ldots,n-1\}$. \qedhere
	\end{enumerate}
\end{proof}

\subsection{Exploiting the Structure of the BRT Component Constraint}
The \ConConstr{} limits the number of upgraded connected components and, as a result, also limits the number of feasible sets of upgraded segments $F \subseteq E$.
This can have an impact on both the number of non-dominated points and the computation time needed to solve the single-objective problems in \Cref{alg:brt-non-dominated}.
For that reason, we further analyze the complexity of the \BRTproblem{} in the context of this component constraint.

First, we consider \BRTtwo{$\star$}{$Z=k$}, where $k$ is fixed and not part of the input, and show that all non-dominated points can be found in polynomial time by an enumeration algorithm, see \Cref{prop:xp}.
This means that the problem is ``slice-wise polynomial'' and, hence, in the complexity class XP \citep{downey2013provable,cygan2015parameterized}.
We remark that \Cref{prop:xp} does not imply that \BRTtwo{$\star$}{$Z\geq 1$} is in FPT, the set of fixed-parameter tractable problems.

\begin{theorem}
	\label{prop:xp}
	\BRTtwo{$\star$}{$Z=k$} can be solved in polynomial time for a fixed $k \in \NN_{\geq 1}$.
\end{theorem}
\begin{proof}
	Let $k\in \NN_{\geq 1}$ be fixed.
	We consider an instance of \BRTtwo{$\star$}{$Z=k$}.
	Each BRT component is uniquely defined by a pair $(s,t) \in V \times V$ with $s < t$ marking the first and last station of the BRT component.
	A feasible set of upgraded segments $F\subseteq E$ can have at most $k$ BRT components determined by $(s_1,t_1),\ldots, (s_k,t_k)$.
	There are at most $n$ possible values for each $s_i$ and $t_i$ for all $i\in \{1,\ldots,k\}$.
	Hence, the number of sets satisfying the \ConConstr{} is in $\mathcal{O}(n^{2k})$.
	For each such set $F$, we can compute the (minimum) investment budget in $\mathcal{O}(\vert E \vert)$ and compute the number of attracted passengers in $\mathcal{O}(\vert D\vert \cdot \vert E\vert)$. 
	Due to \Cref{def:bi-objective-efficient}, a solution $(F,v)$ can only be efficient if $v$ is the minimum investment budget for which $F$ is feasible.
	Therefore, the above procedure gives all potentially efficient solutions of \BRTtwo{$\star$}{$Z=k$}, implying that \BRTtwo{$\star$}{$Z=k$} can be solved in polynomial time for a fixed $k \in \NN_{\geq 1}$.
\end{proof}

Note that the result from \Cref{prop:xp} is especially useful when finding the set of non-dominated points for instances with a low value $k$.

We also consider the case with many BRT components.
Here, \Cref{prop:redundant-component} shows that the \ConConstr{} becomes redundant for values of $Z\geq \left\lceil \frac{\vert E\vert}{2} \right\rceil$.
This lemma thus motivates to consider the case with an arbitrary number of BRT components in more detail. 
Furthermore, \Cref{prop: bound-through-components} shows how we can use \BRT{\Linear}{$Z=\infty$} to obtain bounds on the optimal objective value for the single-objective problem \BRT{\Linear}{$Z=k$} for any fixed $k\in \NN_{\geq 1}$.
Such a bound could, e.g., be used to obtain an approximate Pareto front for instances in which it is hard to solve the single-objective problems in \Cref{alg:brt-non-dominated} to optimality. 

\begin{lemma}
	\label{prop:redundant-component}
	Let an instance of \BRTtwo{$\star$}{$Z = k$} be given.
	The \ConConstr{} is redundant if $Z\geq \left\lceil \frac{\vert E\vert}{2} \right\rceil$.
\end{lemma}

\begin{proof}
	Let an arbitrary subset $F\subseteq E$ be given.
	Then $G[F]$ has the maximum number of connected components if $F$ is a maximum matching in $G$, which would be to take every second segment.
	This yields at most $\left\lceil \frac{\vert E\vert}{2} \right\rceil$ connected components.
	Hence, the number of connected components of $G[F]$ is always less or equal~$Z$ if $Z\geq \left\lceil \frac{\vert E\vert}{2} \right\rceil$.
	In this case, the \ConConstr{} is satisfied.
\end{proof}

\begin{lemma} \label{prop: bound-through-components}
	Let an instance of \BRT{\Linear}{$Z=\infty$} be given.
	Let $f$ and $f_k$ be the optimal objective value of \BRT{\Linear}{$Z=\infty$} with the optimal solution $F$ and of \BRT{\Linear}{$Z=k$} for a fixed $k\in \NN_{\geq 1}$, respectively.
	Let $K$ be the number of BRT components of $F$. 
	Then $f_k=f$ if $k\geq K$, and $f\geq f_k \geq \frac{k}{K} f$ if $k<K$.
\end{lemma}
\begin{proof}
	Dropping the \ConConstr{} is clearly a relaxation, hence, $f\geq f_k$.
	
	If $k\geq K$, then $F$ is still feasible for the restricted problem \BRT{\Linear}{$Z=k$}.
	Therefore, $f_k=f$ in that case.
	
	So let $k<K$, and let $F_1,\ldots,F_K$ be the BRT components of $F$.
	For every $i\in \{1,\ldots,K\}$, we define $r_i \coloneqq \sum_{e\in F_i} \tilde{u}_e$ as the gain in passengers when upgrading the $i$-th BRT component of $F$.
	We assume that they are sorted such that $r_1\geq \ldots \geq r_K\geq 0$.
	Allowing $k$ BRT components means that $F_1\cup \ldots \cup F_k$ is a feasible solution as it has exactly $k$ BRT components.
	This yields that 
	\begin{equation*}
		f_k \geq \sum_{i=1}^k r_i \geq \frac{k}{K} f.
	\end{equation*}
	Here, the last inequality holds because of the following argument:
	Assume that it is not true, which means that ${\sum_{i=1}^k r_i < \frac{k}{K} f}$.
	This implies $\sum_{i=k+1}^K r_i = f-\sum_{i=1}^k r_i > \frac{K-k}{K} f$.
	We then have that $r_k < \frac{1}{K} f$ because we would have $\sum_{i=1}^k r_i \geq k\cdot \frac{1}{K} f$ otherwise, and $r_{k+1}>\frac{1}{K} f$ analogously.
	This is a contradiction to $r_k\geq r_{k+1}$.
\end{proof}

\subsection{Complexity Analysis of the Single-Objective Problem}
While we showed in the previous sections that \SOCtwo{\Linear}{$Z=\infty$} with a special structure and \BRTtwo{$\star$}{$Z=k$} for small values of $k$ can be solved in polynomial time, the time needed to solve the single-objective subproblems \BRT{$\star$}{$Z=\infty$} has a large impact on the running time of \Cref{alg:brt-non-dominated}{} in other cases.
In this section, based on the complexity analysis in \cite{Hoogervorst2022}, we show that the single-objective \BRTproblem{} is related to the well-known knapsack problem and hence NP-hard in general, see \Cref{prop:BRT-Linear-NP} and \Cref{prop:BRT-MinImprov-NP}. 
However, we also identify polynomially solvable cases in \Cref{prop:BRTLinear-poly} and \Cref{prop:BRTMinImprov-poly}. 

\begin{theorem}
	\label{prop:BRT-Linear-NP}
	\BRT{\Linear}{$Z\geq 1$} is NP-hard, even if $Z=\infty$, $\vert M \vert =1$ and $u_e = 1$ for all $e \in E$.
\end{theorem}
\begin{proof}
	For the sake of simplicity, we call the decision version of a problem like its optimization version.
	Given a solution to \BRT{\Linear}{$Z\geq 1$}, we can check in polynomial time whether the budget constraints and the \ConConstr{} are satisfied and a certain value in the objective function is reached.
	
	We reduce (the decision version of) 0-1 knapsack to \BRT{\Linear}{$Z\geq 1$}.
	Let $k$ elements with rewards $r_i\in \NN_{\geq 1}$ and weights $w_i\in \NN_{\geq 1}$ for all $i \in \{1,\ldots,k\}$, a budget $B$ and a bound $S'$ be given.
	We construct an instance of \BRT{\Linear}{$Z\geq 1$} as follows: 
	We set $S\coloneqq S'$, $n\coloneqq k+1$, this means $V\coloneqq\{1,\ldots,k+1\}$, $E\coloneqq \{e_i: i \in\{1,\ldots,k\}\}$, $D\coloneqq \{(i,i+1): i\in \{1,\ldots,k\}\}$, $c_{e_i} \coloneqq w_i$ and $u_{e_i} \coloneqq 1$ for all $i\in\{1,\ldots,k\}$, $M\coloneqq \{1\}$, $b_1\coloneqq 1$, $v \coloneqq B$ and $a_d\coloneqq r_i$ for all $d=(i,i+1)$ with $i\in \{1,\ldots, k\}$ and $Z\coloneqq k$. 
	We show that every feasible solution $F'\subseteq \{1,\ldots,k\}$ of 0-1 knapsack with an objective value of at least $S'$ corresponds to a feasible solution $F\subseteq E$ of \BRT{\Linear}{$Z\geq 1$} with an objective value of at least $S$.
	The solutions $F'$ and $F$ correspond to each other as follows: $i\in F'$ if and only if $e_i \in F$.
	Then the claim holds because
	$\sum_{i\in F'} w_i = \sum_{i\in F'} c_{e_i} = \sum_{e\in F} c_e$
	and 
	\begin{align*}
		\sum_{i\in F'} r_i
		&= \sum_{e_i\in F} a_{(i,i+1)}
		= \sum_{\substack{d=(i,i+1):\\i\in \{1,\ldots,k\}}} \left( \frac{\sum_{e\in F \cap \{e_i\}} 1}{ 1} \cdot a_d \right)\\
		&= \sum_{d\in D} \left( \frac{\sum_{e\in F \cap W_d} u_e}{ \sum_{e\in W_d} u_e} \cdot a_d \right). \qedhere
	\end{align*}
\end{proof}

Consider the mixed-integer programming formulation of \BRT{\Linear}{$Z=\infty$} with a fixed $v\in \RR_{\geq 0}$, where the redundant \ConConstr{}s are dropped:
\begin{maxi*}
	{}{\sum_{e\in E} \tilde{u}_{e} x_e}{\label{ip:linear-reformulation}}{}
	\addConstraint{\sum_{e \in E_m} c_e x_e}{ \leq b_m v}{\quad \text{for all } m \in M}
	\addConstraint{x_e}{ \in \{0,1\}}{\quad \text{for all } e \in E}.
\end{maxi*}

We can see that \BRT{\Linear}{$Z=\infty$} and \SOC{\Linear}{$Z=\infty$} are (multidimensional) 0-1 knapsack problems.
Moreover, because the sets $E_m$, $m\in M$, are disjoint, \BRT{\Linear}{$Z=\infty$} can be decomposed into $\vert M\vert$ independent knapsack problems and hence can be solved in pseudo-polynomial time by dynamic programming.

\begin{theorem}
	\label{prop:BRT-MinImprov-NP}
	\BRT{\MinImprov}{$Z\geq 1$} is NP-hard, even if $Z=\infty$, $M=1$, $u_e=1$ for all $e\in E$ and $L_d=1$ for all $d\in D$.
\end{theorem}
\begin{proof}
	As in the proof of \Cref{prop:BRT-Linear-NP}, \BRT{\MinImprov}{$Z\geq 1$} is in NP.
	
	Further, we apply the same reduction from 0-1 knapsack to \BRT{\MinImprov}{$Z\geq 1$} and additionally choose $L_d \coloneqq 1$ for all $d\in D$.
	It remains to show that the objective value is the same for solutions that correspond to each other.
	We have that 
	\[
	\sum_{\substack{d\in D:\\ L_d\leq \sum\limits_{e\in F\cap W_d}}u_e} a_d
	= \sum_{i\in F'} r_i,
	\]
	because
	\begin{align*}
		\biggl\{ d \in D: L_d & \leq \sum_{e\in F\cap W_d} u_e \biggr\}
		= \biggl\{ (i,i+1): i\in \{1,\ldots, k\} \text{ and } 1\leq \sum_{e\in F\cap \{e_i\}} 1 \biggr\} \\
		&= \{ (i,i+1): i\in \{1,\ldots, k\} \text{ and } e_i \in F\}
		= \{ (i,i+1): i\in F'\}. \qedhere
	\end{align*}
\end{proof}

We conclude the complexity analysis by identifying cases in which the single-objective \BRTproblem{} can be solved in polynomial time.
To this end, we review the consecutive ones property, which is well known in the literature (see, e.g., \cite{Ruf2004,Schoebel2005,Dom2008,Dom2009}).
\Cref{prop:BRTLinear-poly} uses this property to show that \BRT{\Linear}{$Z=\infty$} can be solved in polynomial time in case all segments have unit upgrade costs.
Moreover, \Cref{prop:BRTMinImprov-poly} shows that \BRT{\MinImprov}{$Z=\infty$} can be solved in polynomial time when it holds that all segments have unit upgrade costs and unit improvements, and at the same time only a single segment has to be upgraded to attract the passengers for each OD pair.

\begin{definition}[Consecutive ones property]
	A matrix $A\in \{0,1\}^{k\times l}$ satisfies the consecutive ones property (C1P) on the rows if for all rows $i\in \{1,\ldots,k\}$ it holds: If $A_{i,j}=1$ and $A_{i,j'}=1$ for some $j,j'\in \{1,\ldots, l\}$, $j<j'$, then $A_{i,\bar{j}}=1$ for all $j\leq \bar{j}\leq j'$.
\end{definition}

\begin{lemma}[\citep{Wolsey1999}]
	\label{prop:consecutive-ones}
	If a matrix $A\in \{0,1\}^{k\times l}$ satisfies C1P, then $A$ is totally unimodular.
\end{lemma}

\begin{lemma} \label{prop:BRTLinear-poly}
	\BRT{\Linear}{$Z=\infty$} can be solved in polynomial time if $c_e=1$ for all $e\in E$.
\end{lemma}
\begin{proof}
	Consider
	\begin{maxi}
		{}{\sum_{e\in E} \tilde{u}_{e} x_e}{\label{ip:linear/*/-}}{}
		\addConstraint{\sum_{e \in E_m} c_e x_e}{ \leq b_m v}{\quad \text{for all } m \in M}
		\addConstraint{x_e}{ \in \{0,1\}}{\quad \text{for all }e \in E.}
	\end{maxi}
	We sort the segments and municipalities from one end of the line to the other.
	Let $A\in \RR^{\vert M \vert \times \vert E\vert}$ be the coefficient matrix of the budget constraints, i.e., for all $m\in M$ and $e\in E$, we have $A_{m,e}= 
	1$ if $e\in E_m$, and $A_{m,e}=0$ otherwise.
	Because of the assumption that the municipalities contain only consecutive segments, the matrix $A$ satisfies the consecutive ones property.
	By \Cref{prop:consecutive-ones}, it is totally unimodular and the linear programming relaxation of IP~\eqref{ip:linear/*/-} yields an integer solution.
	Therefore, the problem can be solved in polynomial time~\citep{Wolsey1999}.
\end{proof}

\begin{lemma}
	\label{prop:BRTMinImprov-poly}
	\BRT{\MinImprov}{$Z= \infty$} can be solved in polynomial time if $c_e = 1$, $u_e=1$ for all $e\in E$ and $L_d=1$ for all $d\in D$.
\end{lemma}

\begin{proof}
	We again sort the segments and municipalities from one end of the line to the other.
	The considered special case yields the following simplified formulation:
	\begin{maxi}
		{}{\sum_{d\in D} a_d y_d}{\label{ip:minimprov/*/-}}{}
		\addConstraint{\sum_{e \in E_m} x_e}{ \leq b_m v}{\quad \text{for all } m \in M}
		\addConstraint{\sum_{e \in W_d} \left( -x_e \right) + y_d}{\leq 0}{\quad \text{for all } d \in D}
		\addConstraint{x_e}{ \in \{0,1\}}{\quad \text{for all } e \in E}
		\addConstraint{y_d}{ \in \{0,1\}}{\quad \text{for all } d \in D.}
	\end{maxi}
	The coefficient matrix of the budget constraints and the constraints for the objective is of the form
	$A = \begin{bmatrix}
		A^1 & \textbf{0} \\
		-A^2 & I
	\end{bmatrix},$
	where $I\in \{0,1\}^{\vert D \vert \times \vert D\vert}$ is the unit matrix, $A^1\in \{0,1\}^{\vert M \vert \times \vert E\vert}$ denotes whether a segment belongs to a municipality, and $A^2\in \{0,1\}^{\vert D \vert \times \vert E\vert}$ denotes whether a segment is on the path of an OD pair.
	Formally, we have for all $m\in M$, $d\in D$ and $e\in E$ that \[A^1_{m,e}=\begin{cases}
		1	&\text{if } e\in E_m,\\
		0	&\text{otherwise}
	\end{cases} 
	\text{\quad and \quad}
	A^2_{d,e}=\begin{cases}
		1	&\text{if } e \in W_d,\\
		0	&\text{otherwise.}
	\end{cases}\]
	The matrix $A^1$ has C1P because of the assumption that municipalities contain only consecutive segments, and $A^2$ has C1P because the considered graph is a linear graph.
	As multiplying a row of a matrix by -1 only influences the sign of the determinant of the matrix and its submatrices, the matrix 	$\begin{bmatrix}
		A^1 \\
		-A^2
	\end{bmatrix}$
	is totally unimodular by \Cref{prop:consecutive-ones}.
	This yields that the coefficient matrix $A$, which we obtain by appending a part of a unit matrix to the totally unimodular matrix, is also totally unimodular.
	Therefore, the linear programming relaxation of IP~\eqref{ip:minimprov/*/-} yields an integer solution in this special case, and the problem can be solved in polynomial time~\citep{Wolsey1999}.
\end{proof}

If we consider a global decision maker in addition to the assumptions of \Cref{prop:BRTLinear-poly} and \Cref{prop:BRTMinImprov-poly}, these special cases also satisfy the conditions of \Cref{lem:unit_cost_2-SOC} and thus the complete Pareto front of \SOCtwo{$\star$}{$Z = \infty$} can be constructed in polynomial time.

\section{Numerical Experiments for Artificial Instances}
\label{sec: artificial results}
The Pareto front and the impact of the selected passenger response, the upper bound on the number of BRT components, and the existence of municipalities are at the center of our computational analysis.
These are analyzed in the context of a large library of artificial instances with different interplays between the passenger potential and the upgrade costs. 
Moreover, to investigate the impact of municipalities, we consider different options to split the budget among them.

\subsection{Description of Instances}\label{S.AI:instanceDescr}
All artificial instances consider a line consisting of 25 stations and have the same infrastructure improvements per segment, which are drawn at random. 
The artificial instances differ however in terms of the graph scenario $\alpha = (\alpha_1, \alpha_2)$ defining the costs for upgrading each segment and a demand pattern, respectively, as well as in terms of the budget scenario $\beta$ defining the budget split among five municipalities.
The considered values for these parameters are given in \Cref{tab:param}.

The cost pattern $\alpha_1$ varies between uniform costs per segment (UNIT), a pattern with higher costs in the center of the line (MIDDLE), and an ENDS pattern to contrast with the others by having the segments with the highest upgrade costs at the ends of the line.
The MIDDLE cost pattern could result, e.g., from upgrades in the inner city being more complicated, while ENDS could represent outside-city upgrades involving highway lanes that are very expensive to upgrade.
The cost patterns together with the (fixed) infrastructure improvements are depicted in \Cref{fig:costs} in the Appendix.

\begin{figure}[h!t]
	\centering
	\begin{subfigure}[b]{0.475\textwidth}
		\centering
		\includegraphics[width=0.95\textwidth]{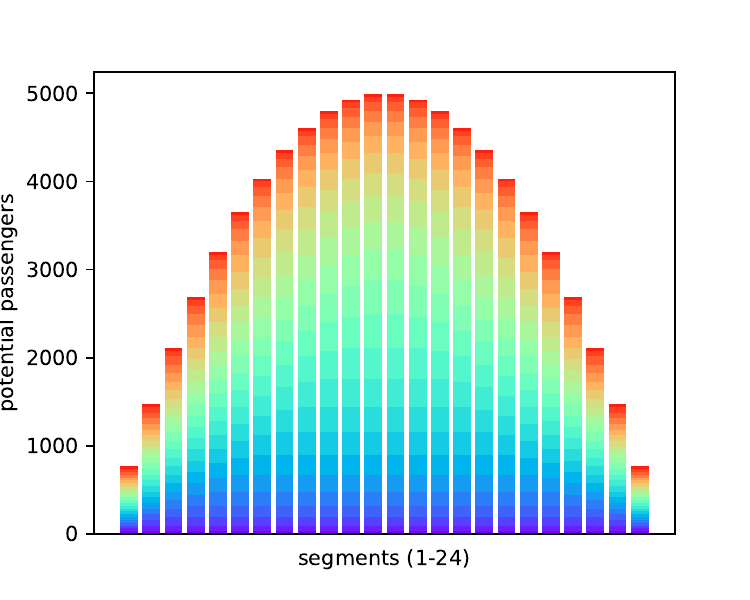}
		\caption{\EvenDemand}
		\label{fig:demand_even}
	\end{subfigure}
	\begin{subfigure}[b]{0.475\textwidth}
		\centering
		\includegraphics[width=0.95\textwidth]{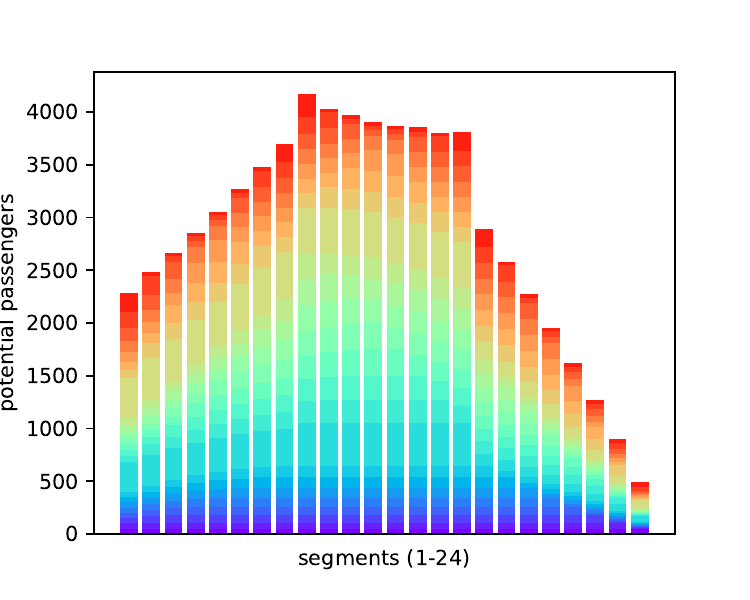}
		\caption{\LargeStationsDemand}
		\label{fig:demand_center}
	\end{subfigure}
	\begin{subfigure}[b]{0.475\textwidth}
		\centering
		\includegraphics[width=0.95\textwidth]{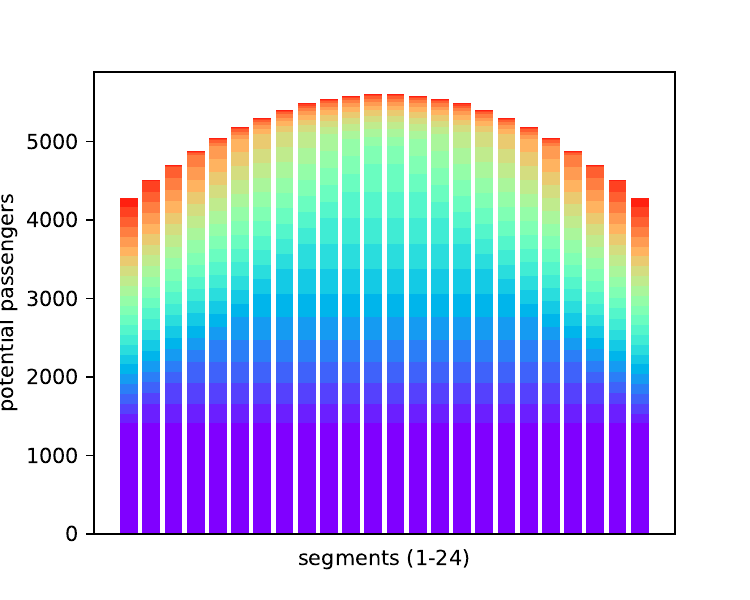}
		\caption{\EndStationsDemand}
		\label{fig:demand_end}
	\end{subfigure}
	\begin{subfigure}[b]{0.475\textwidth}
		\centering
		\includegraphics[width=0.95\textwidth]{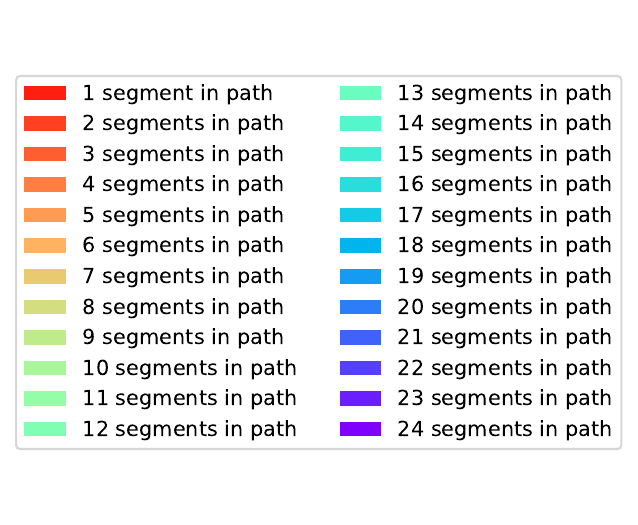}
		\caption{Legend}
	\end{subfigure}
	\caption{Load profiles of all demand patterns for the artificial instances. The horizontal axis contains each of the 24 segments of the BRT line. Each bar represents the number of potential passengers using the respective segment, and the coloring of a bar depicts the total travel distance of passengers using that segment.}
	\label{fig:demand}
\end{figure}

Also the volume and the distribution of the potential passengers of the line impact the solutions. 
The three different load profiles resulting from the three demand patterns $\alpha_2$ (\EvenDemand{}, \LargeStationsDemand{}, \EndStationsDemand{}) are depicted in \Cref{fig:demand}, with the height of the bar indicating the load per segment, and the colored shading indicating the length of the boarded passengers' paths. 
Thus, \LargeStationsDemand{} results typically in shorter path lengths than \EndStationsDemand{} and \EvenDemand{},
whereas \EndStationsDemand{} has especially many passengers traveling from one end station to the other one.
Moreover, \EvenDemand{} has fewer passengers traveling around the terminals of the line than the other two.
For additional information regarding the locations of stations with high demand along the line and the travel distances of passengers in these demand patterns, see \Cref{fig:graphs} and \Cref{fig:hist-distances} in the Appendix.

The budget split $\beta$ describes the distribution of the total available investment budget among the municipalities. 
We consider a distribution according to equal shares (EQUAL), proportional to the costs required for upgrading all segments in a municipality (COST), and according to the total passenger volume that flows in and out of the stations belonging to the municipality (PASS).

\begin{table}[h!t]
	\caption{Parameters for generating artificial instances.}
	\label{tab:param}
	\centering
	\begin{tabular}{p{0.025\textwidth} @{\hskip 0.75\tabcolsep}p{0.19\textwidth}@{\hskip 0.75\tabcolsep}p{0.105\textwidth}p{0.61\textwidth}}
		\toprule
		\multicolumn{2}{l}{Parameter} & Value & Explanation\\ 
		\midrule
		\multirow{3}{*}{$\alpha_1$} & \multirow{3}{*}{cost pattern} &  UNI & unit costs of $c_e = 1$ for all $e \in E$ \\
		&  &  MIDDLE & more expensive towards the middle of the line\\
		&  &  ENDS & more expensive towards the end stations of the line\\
		\addlinespace
		\multirow{4}{*}{$\alpha_2$} & \multirow{4}{*}{demand pattern} &  \EvenDemand{}  & same passenger potential for each OD pair \\
		&  &   \LargeStationsDemand{} & centered around large stations, passengers distributed according to the gravity model \citep{Rodrigue2020} \\
		&  &  \EndStationsDemand{} & high passenger potential between end stations of the line \\
		\addlinespace
		\midrule 
		\addlinespace
		\multirow{6}{*}{$\beta$} & \multirow{6}{*}{budget split} & EQUAL & budget distributed equally among municipalities, i.e., equal budget shares $b_m$ \\	
		& & COST & budget shares $b_m$ proportional to the costs of the segments in municipality $m$\\	
		& & PASS & budget shares $b_m$ proportional to the number of potential passengers entering or exiting in municipality $m$\\		
		\bottomrule
	\end{tabular}
\end{table}

\subsection{Computational Study Design}
\label{sec: instances}

In our computational study, we compute the Pareto fronts for all instances, focusing on the part where the investment budget does not exceed the investment costs for upgrading all segments, i.e., $v \leq c(E)$.
As a consequence, dependent on the budget split, not all municipalities may have sufficient budget to upgrade all of their segments.
In addition to varying characteristics of the instances (\Cref{S.AI:instanceDescr}), we consider the different problem variants \BRTfullTwo{$\lambda_1$}{$\lambda_2$}{$\lambda_3$} from \Cref{sec:problem-variants}:

\paragraph{Passenger Response}
All instances are evaluated for the objectives \Linear{} and \MinImprov{}. 
For \MinImprov{}, we require that a minimum of 75\% of the potential infrastructure improvements is achieved through upgrades before the passengers corresponding to that OD pair are attracted, i.e., 
\begin{equation*}
	L_d \coloneqq \left\lfloor 0.75 \cdot \sum_{e \in W_d} u_e \right\rfloor
\end{equation*}
for all OD pairs $d \in D$.

\paragraph{Number of BRT components}
We consider upper bounds on the number of BRT components $Z\in \{1,2,3,\infty\}$.
Our experiments show that the difference between $Z=3$ and $Z=\infty$ is generally small, and therefore including more options for~$Z$ than $\{1,2,3,\infty\}$ would not lead to further insights in our setting. 

\paragraph{Municipalities}
In order to determine the impact of the separate municipality budgets, each instance is evaluated both in the context of a global decision maker with a single budget ($|M| = 1$) as well as in the original context where each municipality has its own budget constraint ($|M| = 5$).
In the former setting, the global decision maker can spend the whole investment budget $v$, i.e., there is a single municipality with $b_1 = 1$, while in the latter setting, the investment budget is distributed among the municipalities according to the budget split~$\beta$ (\Cref{tab:param}).

The data of the artificial instances together with the applied upper bound on the number of BRT components and the municipality scenario, i.e., a total of $3^2\cdot 4^2=144$ settings that are evaluated regarding both passenger responses, is available at \url{https://doi.org/10.11583/DTU.23653893}.

\subsection{Computation Time}
\label{sec:computation-time}
All instances are solved by means of the commercial solver CPLEX 22.1 on a computer with an Intel Xeon Gold 6126 processor, using 12 CPU cores and a total of 24 GB of internal memory.
The corresponding running time for computing the non-dominated points for the artificial instances is shown in \Cref{tab:runtime_25}.
Here, we give the average time to find the Pareto front, the average number of points on the Pareto front and the average time for obtaining a single non-dominated point for each passenger response, each cost pattern $\alpha_1$ and both municipality scenarios $\vert M \vert \in \{1,5\}$.
Note that the reported values are averaged over all three demand patterns $\alpha_2$ and over the considered upper bound on the number of BRT components $Z \in \{1,2,3,\infty\}$.
Additionally, the results are also averaged over the different budget splits $\beta$ for the setting with municipalities ($|M| =5$).

\begin{table}[tb]
	\caption{Running time in seconds, number of obtained Pareto points and running time per Pareto point for problem variants \SOCtwo{$\star$}{$Z \geq 1$} and \BRTfullTwo{$\star$}{$Z \geq 1$}{$\vert M \vert =5$}. The results have been averaged over artificial instances sharing the same cost pattern $\alpha_1$.}
	\label{tab:runtime_25}
	\centering
	\resizebox{\textwidth}{!}{
		\begin{tabular}{llrrrrrr}
			\toprule
			&& \multicolumn{3}{c}{\SOCtwo{$\star$}{$Z \geq 1$}} & \multicolumn{3}{c}{\BRTfullTwo{$\star$}{$Z \geq 1$}{$\vert M \vert =5$}}\\
			\cmidrule(rl){3-5} \cmidrule(rl){6-8}
			Objective & $\alpha_1$  & all points & \# points & per point & all points & \# points & per point \\
			\midrule
			\Linear{}    & UNIT   & 0.17               & 25.00              & 0.007             & 0.09               & 7.47               & 0.012              \\
			\Linear{}    & MIDDLE & 3.53               & 179.83            & 0.020              & 0.52               & 37.53              & 0.014              \\
			\Linear{}    & ENDS   & 1.95               & 107.00             & 0.016             & 0.39               & 28.25              & 0.014              \\
			\midrule
			\MinImprov{} & UNIT   & 14.96              & 25.00              & 0.599             & 1.12               & 7.42               & 0.152              \\
			\MinImprov{} & MIDDLE & 768.00              & 105.42            & 5.098             & 3.93               & 26.08              & 0.139              \\
			\MinImprov{} & ENDS   & 68.65              & 77.17             & 0.639             & 3.02               & 25.00               & 0.117   \\            
			\bottomrule
	\end{tabular}}
\end{table}

The results in \Cref{tab:runtime_25} show that the Pareto fronts can overall be computed quickly, especially for the \Linear{} passenger response.
Moreover, it can be seen that the introduction of separate municipality budgets ($|M| =5$) consistently leads to a lower running time and fewer points on the Pareto front than the consideration of a global decision maker ($|M| = 1$).
This is likely a result of the smaller solution space with separate municipality budgets, where fewer combinations of items fit within the individual municipality budgets.
The longest running times can be observed for the MIDDLE cost pattern in combination with the \MinImprov{} passenger response, where especially the long running time for the setting with a global decision maker ($|M|=1$) stands out.
This might be explained by the middle segments often having the highest passenger load as well as being the most expensive to upgrade when considering this cost pattern.
Looking at the number of points on the Pareto front, it can be observed that the number of non-dominated points is significantly lower than the theoretical bound determined in \Cref{lem:non-dominated-bound}.
Moreover, in line with \Cref{lem:unit_cost_2-SOC}, it can be seen that the number of non-dominated points is often significantly lower for the UNIT cost pattern than for the other cost patterns.

\begin{figure}[ht]
	\centering
	\begin{subfigure}[b]{0.4\textwidth}
		\centering
		\includegraphics[width=0.95\textwidth]{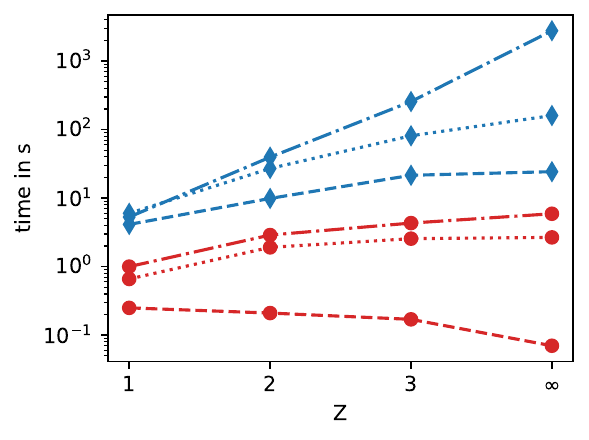}
		\caption{$\vert M \vert = 1$}
		\label{fig:runtime-total-time-1}
	\end{subfigure}
	\begin{subfigure}[b]{0.4\textwidth}
		\centering
		\includegraphics[width=0.95\textwidth]{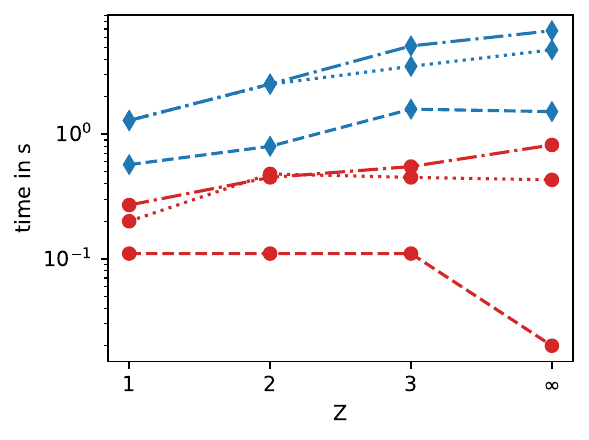}
		\caption{$\vert M \vert = 5$}
		\label{fig:runtime-total-time-5}
	\end{subfigure}
	\begin{subfigure}[b]{0.18\textwidth}
		\centering
		\raisebox{15mm}{\includegraphics[width=0.95\textwidth]{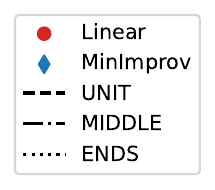}}
		\caption{Legend}
	\end{subfigure}
	\caption{Running time of \BRTfullTwo{$\star$}{$\vert M \vert =1$}{$Z\geq1$} and \BRTfullTwo{$\star$}{$\vert M\vert=5$}{$Z\geq 1$} with a logarithmic scale on the vertical axis. The values are averaged over all demand patterns and budget splits (if applicable).}
	\label{fig:runtime-total-time}
\end{figure}

\begin{figure}[ht]
	\centering
	\begin{subfigure}[b]{0.4\textwidth}
		\centering
		\includegraphics[width=0.95\textwidth]{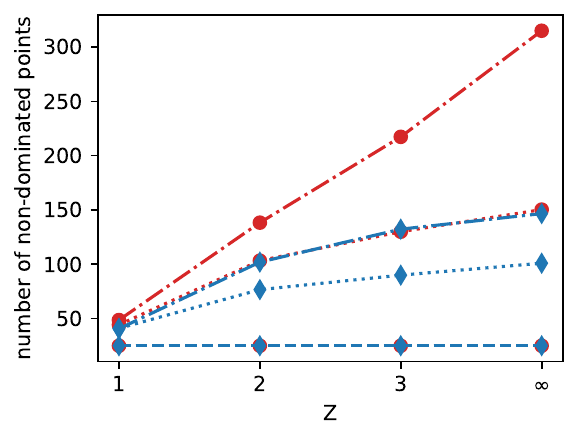}
		\caption{$\vert M \vert = 1$}
		\label{fig:runtime-number-points-1}
	\end{subfigure}
	\begin{subfigure}[b]{0.4\textwidth}
		\centering
		\includegraphics[width=0.95\textwidth]{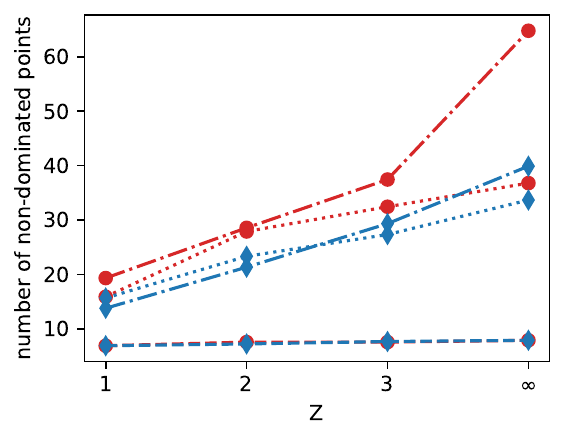}
		\caption{$\vert M \vert = 5$}
		\label{fig:runtime-number-points-5}
	\end{subfigure}
	\begin{subfigure}[b]{0.18\textwidth}
		\centering
		\raisebox{15mm}{\includegraphics[width=0.95\textwidth]{runtime_legend.pdf}}
		\caption{Legend}
	\end{subfigure}
	\caption{Number of non-dominated points of \BRTfullTwo{$\star$}{$\vert M \vert =1$}{$Z\geq1$} and \BRTfullTwo{$\star$}{$\vert M\vert=5$}{$Z\geq 1$}. The values are averaged over all demand patterns and budget splits (if applicable).}
	\label{fig:runtime-number-points}
\end{figure}

In addition, \Cref{fig:runtime-total-time} shows the average running time dependent on the upper limit $Z$ on the number of BRT components.
We can see that the average running time increases with $Z$, which is mainly because of the growing number of non-dominated points (see \Cref{fig:runtime-number-points}), except for the case of \Linear{} together with the cost pattern UNIT, in which it decreases.
\Cref{prop:BRTLinear-poly} showed that this setting together with $Z=\infty$ is a polynomial time special case.

\subsection{Analysis Of Pareto Fronts}

\begin{figure}[htbp]
	\centering
	
	\makebox{\raisebox{80pt}{\rotatebox[origin=c]{90}{$\alpha_1=\textrm{UNIT}$}}}%
	\;
	\begin{subfigure}[b]{0.31\textwidth}
		\centering
		\includegraphics[width=\textwidth]{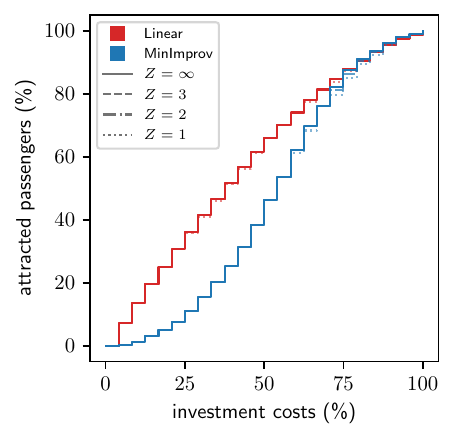}
	\end{subfigure}
	\begin{subfigure}[b]{0.31\textwidth}
		\centering
		\includegraphics[width=\textwidth]{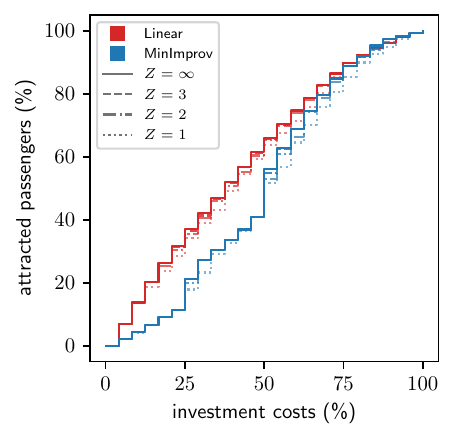}
	\end{subfigure}
	\begin{subfigure}[b]{0.31\textwidth}
		\centering
		\includegraphics[width=\textwidth]{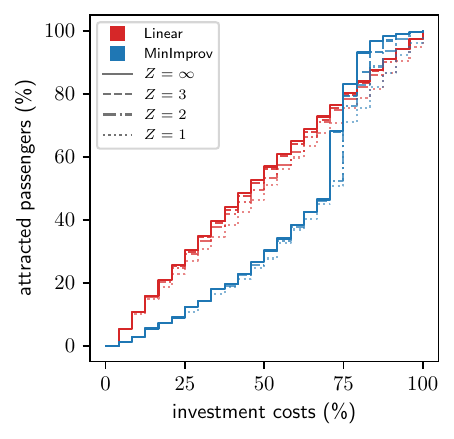}
	\end{subfigure}
	
	\makebox{\raisebox{85pt}{\rotatebox[origin=c]{90}{$\alpha_1=\textrm{MIDDLE}$}}}%
	\;
	\begin{subfigure}[b]{0.31\textwidth}
		\centering
		\includegraphics[width=\textwidth]{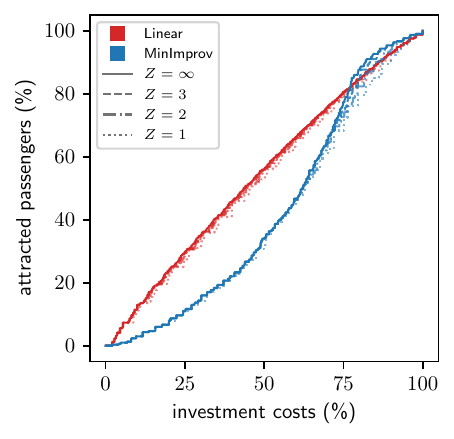}
	\end{subfigure}
	\begin{subfigure}[b]{0.31\textwidth}
		\centering
		\includegraphics[width=\textwidth]{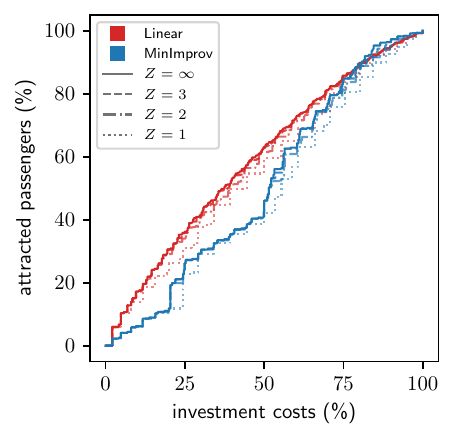}
	\end{subfigure}
	\begin{subfigure}[b]{0.31\textwidth}
		\centering
		\includegraphics[width=\textwidth]{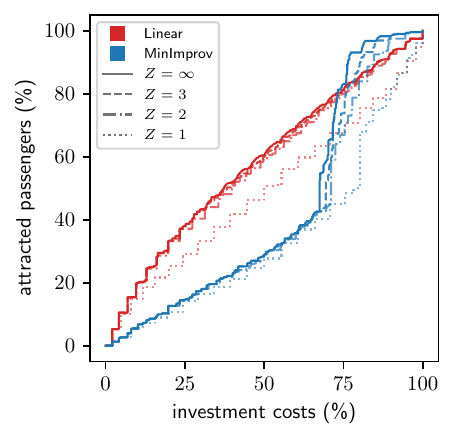}
	\end{subfigure}
	
	\makebox{\raisebox{85pt}{\rotatebox[origin=c]{90}{$\alpha_1=\textrm{ENDS}$}}}%
	\;
	\begin{subfigure}[b]{0.31\textwidth}
		\centering
		\includegraphics[width=\textwidth]{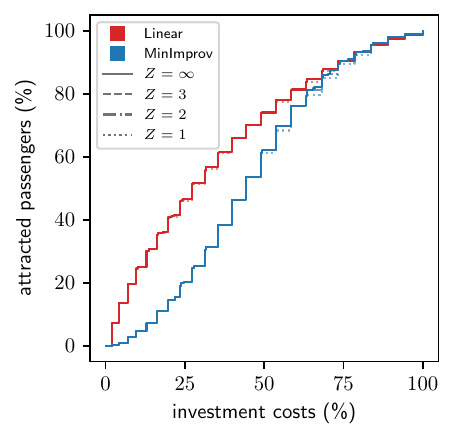}
		\captionsetup{justification=centering}
		\caption*{$\alpha_2=\textrm{\EvenDemand}$}
	\end{subfigure}
	\begin{subfigure}[b]{0.31\textwidth}
		\centering
		\includegraphics[width=\textwidth]{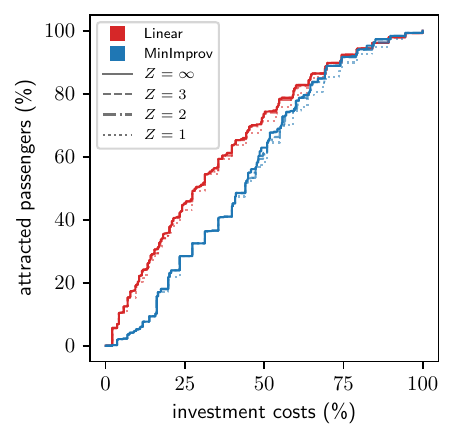}
		\captionsetup{justification=centering}
		\caption*{$\alpha_2=\textrm{\LargeStationsDemand}$}
	\end{subfigure}
	\begin{subfigure}[b]{0.31\textwidth}
		\centering
		\includegraphics[width=\textwidth]{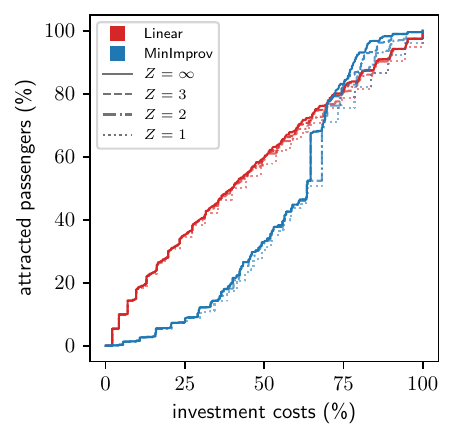}
		\captionsetup{justification=centering}
		\caption*{$\alpha_2=\textrm{\EndStationsDemand}$}
	\end{subfigure}
	
	\caption{Evaluation of the non-dominated points of \SOCtwo{\Linear}{$Z\geq1$} (red) and \SOCtwo{\MinImprov}{$Z\geq 1$} (blue) for artificial instances representing all choices for parameters $\alpha_1, \alpha_2$ and $Z$. Both attracted passengers and investment costs are given as percentage of the total number of potential passengers and costs for upgrading all segments, respectively.}
	\label{fig:pareto_compare_obj}
\end{figure}

In this section, we analyze the influence of the passenger response, the number of BRT components, the demand pattern and the municipalities on the Pareto front.
As described in \Cref{sec:eval_investment}, we compute the efficient solutions and the Pareto fronts with respect to the investment budget, but we evaluate the results with respect to the investment costs.
Therefore, the following figures show the investment costs on the horizontal axis and the newly attracted passengers on the vertical axis.
Both are given as percentage of the total number of potential passengers and costs for upgrading all segments, respectively.
\Cref{fig:pareto_compare_obj} shows the evaluation for a global decision maker ($\vert M \vert =1$).
The red plots represent the \Linear{} passenger response and the blue plots represent the \MinImprov{} passenger response, while the line style represents the number of allowed BRT components $Z$.
All graphs in a row share the same cost pattern $\alpha_1$, and all graphs in a column share the same demand pattern $\alpha_2$.

\paragraph*{Influence of the Passenger Response}
In general, the non-linear objective \MinImprov{} leads to solutions with less passengers per investment budget than \Linear{}, with the exception of high level investments of at least around 75\% of the total budget or more. 
The cut-off point at 75\% correlates with the minimum improvement $L_d$ of 75\% required within \MinImprov{}, for $(\alpha_1, \alpha_2) = (\textrm{ENDS}, \textrm{\EndStationsDemand{}})$ the cut-off point is a bit lower.
Furthermore, the shape of the curve is typically more convex for \MinImprov{}, in which the return on investment is generally increasing and only starts to reduce much later than for the \Linear{} passenger response.
\Linear{} rather shows a higher return on investments at the lower investment levels.
This can be explained by looking at the passenger responses.
For the \Linear{} passenger response, passengers of all OD pairs that are affected by upgrades are attracted in proportion to the realized infrastructure improvements.
For the \MinImprov{} passenger response, however, mainly passengers of affected OD pairs with a short travel distance are attracted at a low investment budget level.
Only at higher investment budget levels, when sufficiently many segments can be upgraded, long-distance travelers are also attracted.
As the demand pattern \EndStationsDemand{} has around 14\% of all passengers traveling along all 24 segments and the aggregated demand over the other OD pairs decreases only slowly in the path length (see \Cref{fig:hist-distances} in the Appendix), the convexity effect is most strongly pronounced for this demand pattern.
In comparison, the convexity effect is a bit less pronounced for \EvenDemand{} and only weakly present for \LargeStationsDemand{}, for which the demand of long-distance journeys is generally decreasing.
The gap between the two passenger responses is generally smaller for high investment budgets. 
These results indicate that the passenger response has a strong impact on the trade-off between attracted passengers and investments.
Investigating the passenger behavior as part of BRT feasibility studies would thus be important to determine an appropriate investment level.

\paragraph*{Influence of the Number of BRT Components}
The impact of the upper bound on the number of BRT components $Z$ diminishes quickly with size, where the numbers of attracted passengers and the investment costs of non-dominated solutions for $Z=3$ and $Z=\infty$ are almost identical.
Again, the impact of $Z$ is higher for larger investment budgets and also more prevailing for the \MinImprov{} passenger response. 
Additionally, the MIDDLE cost pattern (center row) and the \EndStationsDemand{} demand pattern (right column) show a large impact of the \ConConstr{}. 
In general, we see that restricting the number of BRT components to a fixed $Z>1$ comes at small costs, while it could lead to lines that may be considered of higher quality from a passenger perspective. 
Finally, from a computational perspective fixing $Z$ can reduce the computational complexity, as shown earlier in \Cref{prop:xp} and in \Cref{sec:computation-time}.

\begin{figure}[ht]
	\centering
	\begin{subfigure}[b]{0.475\textwidth}
		\centering
		\includegraphics[width=0.95\textwidth]{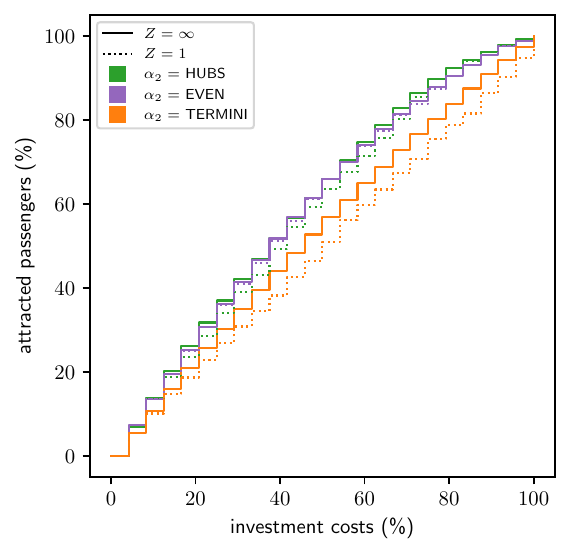}
		\caption{\Linear{}}
		\label{fig:pareto_25_unit:Lin}
	\end{subfigure}
	\begin{subfigure}[b]{0.475\textwidth}
		\centering
		\includegraphics[width=0.95\textwidth]{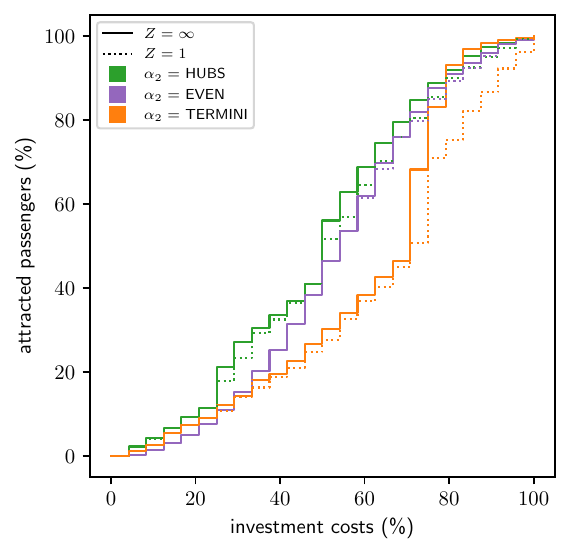}
		\caption{\MinImprov{}}
		\label{fig:pareto_25_unit:Min}
	\end{subfigure}
	\caption{Evaluation of the non-dominated points of $\SOCtwo{\star}{Z\geq 1}$ for artificial instances with cost pattern $\alpha_1 = \textrm{UNIT}$ and $Z \in \{1,\infty\}$ and all choices for the demand pattern~$\alpha_2$. Both attracted passengers and investment costs are given as percentage of the total number of potential passengers and costs for upgrading all segments, respectively.}
	\label{fig:pareto_25_unit}
\end{figure}

\paragraph*{Influence of the Demand Pattern}
\Cref{fig:pareto_25_unit} depicts the effect of the demand pattern on the sets of non-dominated points for $\SOCtwo{\star}{Z=\infty}$ (solid lines) and $\SOCtwo{\star}{Z=1}$ (dotted lines) for the cost pattern UNIT (the results for the other cost patterns are similar and can be seen in \Cref{fig:pareto_25_cost_pattern} in the Appendix).
For the passenger response \Linear{}, the demand patterns behave similarly.
The only thing that stands out is that \EndStationsDemand{} leads to slightly fewer attracted passengers compared to \LargeStationsDemand{} and \EvenDemand{}.
This can also be seen for \MinImprov{}.
In addition, we see a large jump in attracted passengers for \MinImprov{} with demand pattern \EndStationsDemand{} when around 75\% of the budget is invested.
This is due to the relatively high number of passengers that travel along all 24 segments (about 14\% of all passengers) because realizing 75\% of the potential improvement suffices to attract all those passengers according to the definition of \MinImprov{}. 
The influence of restricting the number of connected components to $Z=1$ is especially pronounced for the demand pattern \EndStationsDemand{}. 
Here again, the high number of passengers using all 24 segments is affected most by restricting the set of upgraded segments.

\begin{figure}[htbp]
	\centering
	
	\makebox{\raisebox{80pt}{\rotatebox[origin=c]{90}{$\alpha_1=\textrm{UNIT}$}}}%
	\;
	\begin{subfigure}[b]{0.31\textwidth}
		\centering
		\includegraphics[width=\textwidth]{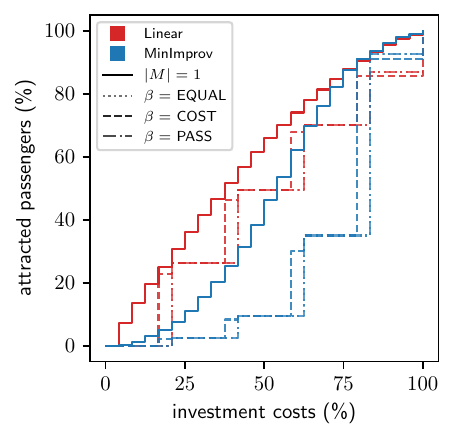}
	\end{subfigure}
	\begin{subfigure}[b]{0.31\textwidth}
		\centering
		\includegraphics[width=\textwidth]{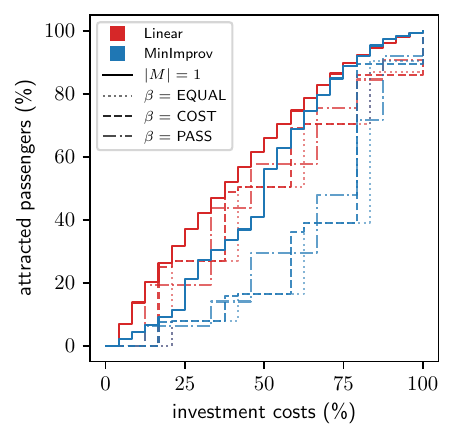}
	\end{subfigure}
	\begin{subfigure}[b]{0.31\textwidth}
		\centering
		\includegraphics[width=\textwidth]{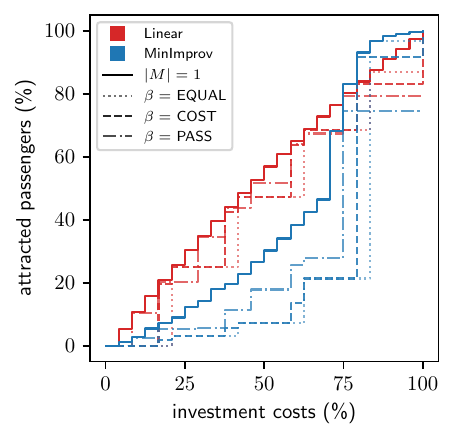}
	\end{subfigure}

	\makebox{\raisebox{85pt}{\rotatebox[origin=c]{90}{$\alpha_1=\textrm{MIDDLE}$}}}%
	\;
	\begin{subfigure}[b]{0.31\textwidth}
		\centering
		\includegraphics[width=\textwidth]{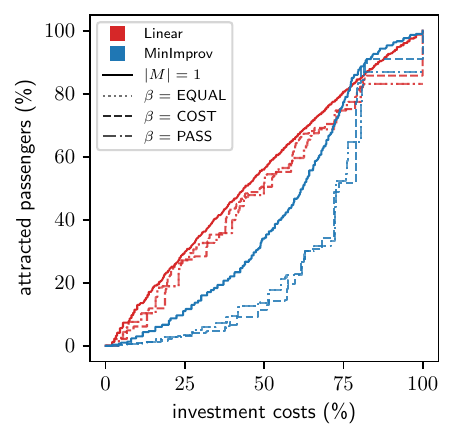}
	\end{subfigure}
	\begin{subfigure}[b]{0.31\textwidth}
		\centering
		\includegraphics[width=\textwidth]{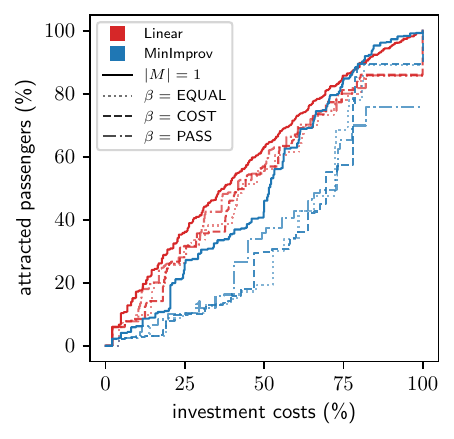}
	\end{subfigure}
	\begin{subfigure}[b]{0.31\textwidth}
		\centering
		\includegraphics[width=\textwidth]{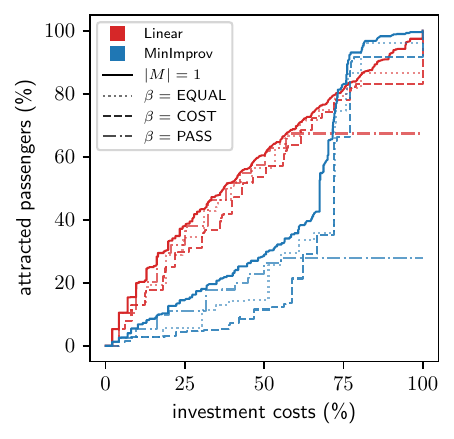}
	\end{subfigure}

	\makebox{\raisebox{85pt}{\rotatebox[origin=c]{90}{$\alpha_1=\textrm{ENDS}$}}}%
	\;
	\begin{subfigure}[b]{0.31\textwidth}
		\centering
		\includegraphics[width=\textwidth]{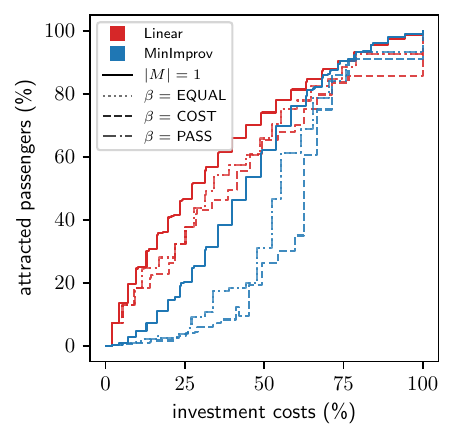}
		\captionsetup{justification=centering}
		\caption*{$\alpha_2=\textrm{\EvenDemand}$}
	\end{subfigure}
	\begin{subfigure}[b]{0.31\textwidth}
		\centering
		\includegraphics[width=\textwidth]{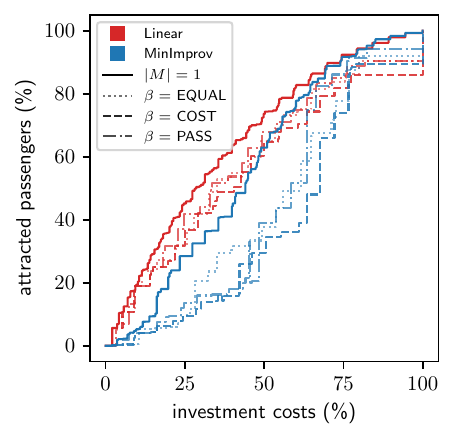}
		\caption*{$\alpha_2=\textrm{\LargeStationsDemand}$}
	\end{subfigure}
	\begin{subfigure}[b]{0.31\textwidth}
		\centering
		\includegraphics[width=\textwidth]{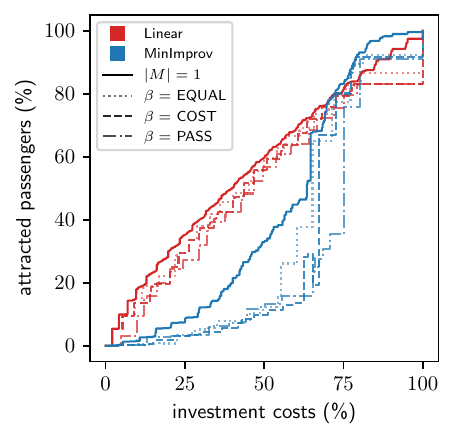}
		\captionsetup{justification=centering}
		\caption*{$\alpha_2=\textrm{\EndStationsDemand}$}
	\end{subfigure}
	
	\caption{Evaluation of the non-dominated points of $\BRTtwo{\Linear}{Z=\infty}$ (red) and $\BRTtwo{\MinImprov}{Z=\infty}$ (blue) for artificial instances with $\vert M \vert \in \{1,5\}$ and all choices for parameters $\alpha_1, \alpha_2$ and $\beta$. Both attracted passengers and investment costs are given as percentage of the total number of potential passengers and costs for upgrading all segments, respectively.}
	\label{fig:pareto_compare_objMuni}
\end{figure}

\paragraph*{Influence of Municipalities}
The impact of the distribution of the budget among the municipalities is depicted in \Cref{fig:pareto_compare_objMuni}, which is similar in set-up to \Cref{fig:pareto_compare_obj} with the difference that the line styles now represent different budget splits among the municipalities, with the solid line representing the case of a global decision maker. 
Moreover, all results in \Cref{fig:pareto_compare_objMuni} are obtained without the \ConConstr{} ($Z=\infty$).

The introduction of municipalities generally leads to lower numbers of attracted passengers. 
Because of the distribution of the investment budget among the municipalities, compared to the case of a global decision maker, only a smaller share can be invested and not always in the segments that would attract the most passengers.
Moreover, considering several municipalities intensifies the findings of the case with a global decision maker: \MinImprov{} requires higher investments for the same number of passengers until around 75\% of investments and is characterized by a return on investment that follows a more convex shape compared to the \Linear{} passenger response, with the same explanations as for the case of a global decision maker.
The impact of the chosen budget split among the municipalities is typically higher for \MinImprov{} as well. 
For budget splits other than COST, the full upgrade may not be achievable even at an investment budget equal to 100\% of the total upgrade costs because individual municipalities may not have enough money available to upgrade all segments belonging to them.
This is specifically visible for $\alpha_1=\text{MIDDLE}$, $\alpha_2=\text{\EndStationsDemand}$, $\beta= \textrm{PASS}$, in this case, even only about 30\% of the passenger potential are attracted at 100\% investment budget, showing that the investment costs stagnate at 55\% because of the interplay between budget split and the demand pattern.
The results indicate that, especially in case of a non-linear relationship between BRT upgrades and attracted passengers, establishing a framework for collaboration and co-commitment has a large influence on the number of attracted passengers and thereby on the line potential. 

\section{Greater Copenhagen Case Study}
\label{sec: case study results}
We now focus on the case study for the planned BRT line in Greater Copenhagen.
We first describe the case study and the corresponding instances.
Afterwards, we analyze the Pareto plots that are obtained for these instances.

\subsection{Description of the Instances}
\label{sec:Copenhagen-case}
\begin{figure}[!ht]
	\centering
	\includegraphics[width=0.65\textwidth]{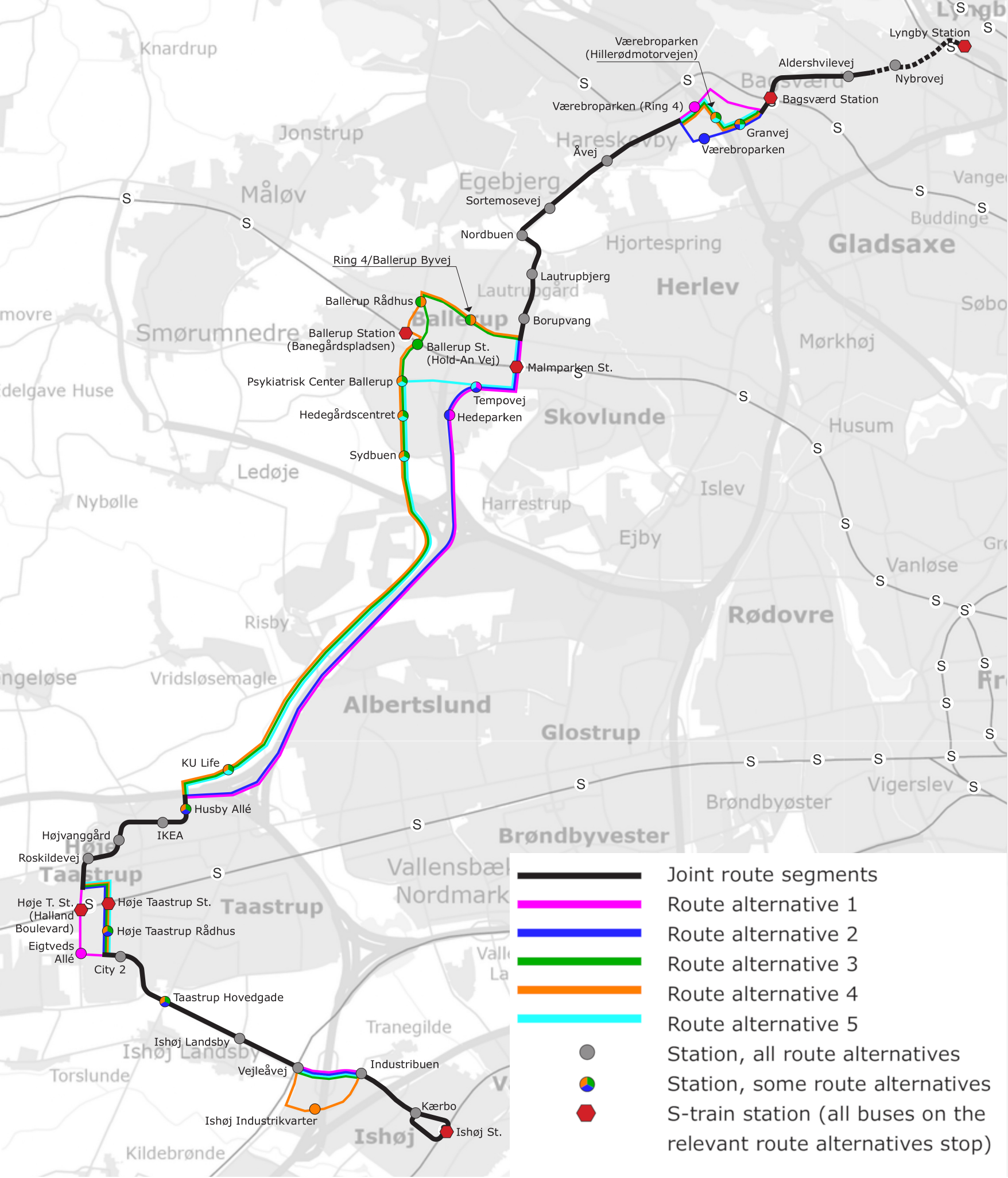}
	\caption{Route alternatives for a new BRT line in Greater Copenhagen. Adapted from \cite{vejdirektoratet2022brt}.}
	\label{fig:brt}
\end{figure}

Currently, the Capital Region in Denmark (a regional government) is planning to build a set of BRT lines within Copenhagen and the urban area surrounding it, i.e., Greater Copenhagen.
One of these new BRT lines will run foremost along the route of the bus line 400S, which is currently a traditional mixed traffic line. 
A pre-assessment study was conducted for the BRT line that calculated the expected costs, travel durations and number of passengers per station for five different route alternatives \citep{vejdirektoratet2022brt}.
These five route alternatives are shown in \Cref{fig:brt}.

All five route alternatives run through a total of eight municipalities.
These municipalities have authority over their local investments in the BRT line, and investments in local infrastructure would not be possible without their involvement. 
Next to the municipalities, also the Capital Region, the central Danish government, and Movia are involved in the planning process.
Movia is a public transport agency funded by the collective of municipalities in the Capital Region, which highlights the willingness of the municipalities to work together to find socially optimal solutions for public transport in the region.
Moreover, due to the expertise available within the agency, Movia overall takes a leading role in the design of the new BRT line and thereby provides suggestions that then need to be approved by the municipalities.
This process can be iterative: Municipalities discuss solutions and revise their budget levels, followed by new suggestions from Movia.
Hereby, the proposed model can advise Movia on how sensitively the number of newly attracted passengers reacts to a reduction in the investment budget. 
Thereby, it can illustrate the importance of achieving a high upgrade level and may aid the transport agency in selecting its strategy.

We use the data from the pre-assessment study to derive instances for the \BRTproblem{} for each of the five route alternatives. 
These instances contain between 24 and 32 stations, depending on the route. 
The current plan includes connecting the BRT line via Nybrovej to Lyngby station, even though the responsible municipality has indicated it is not willing to invest in upgrading the infrastructure on their segments. 
Therefore, our case includes two segments that cannot be upgraded.
The remaining seven municipalities are willing to partake in the BRT project.
The upgrade costs per segment are derived from the required infrastructure investments for the line provided in the pre-assessment. 
Moreover, the potential benefit of upgrading a segment is defined by the difference between expected travel time of the current mixed traffic line and the new expected travel times of the BRT line as defined in the pre-assessment.
The upgrade costs and infrastructure improvements are depicted in \Cref{fig:costs-movia} in the Appendix.

\begin{figure}[tbp]
	\centering
	\begin{subfigure}{0.32 \textwidth}
		\includegraphics[width=\textwidth]{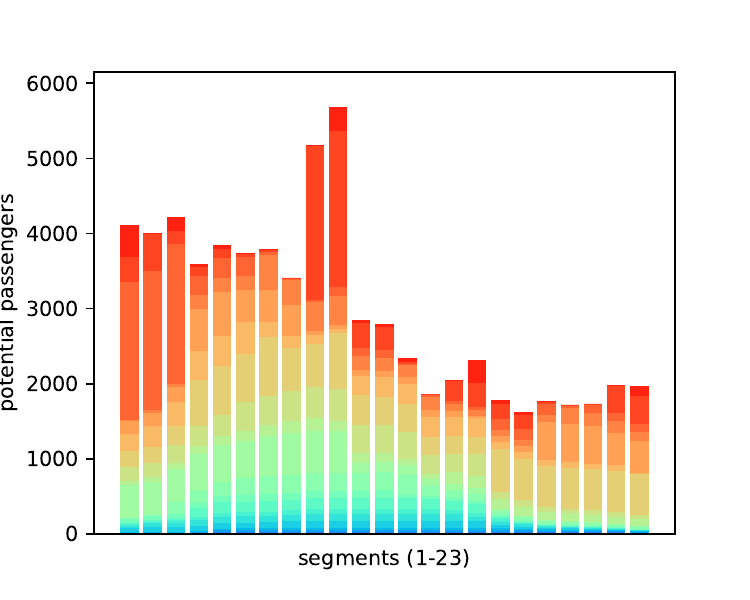}
		\caption{Alternative 1}
	\end{subfigure}
	\begin{subfigure}{0.32 \textwidth}
		\centering
		\includegraphics[width=\textwidth]{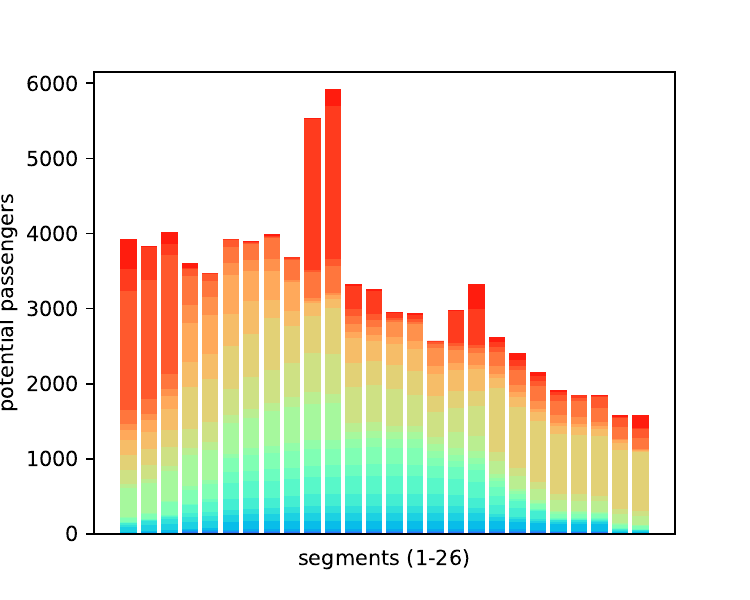}
		\caption{Alternative 2}
	\end{subfigure}
	\begin{subfigure}{0.32 \textwidth}
		\centering
		\includegraphics[width=\textwidth]{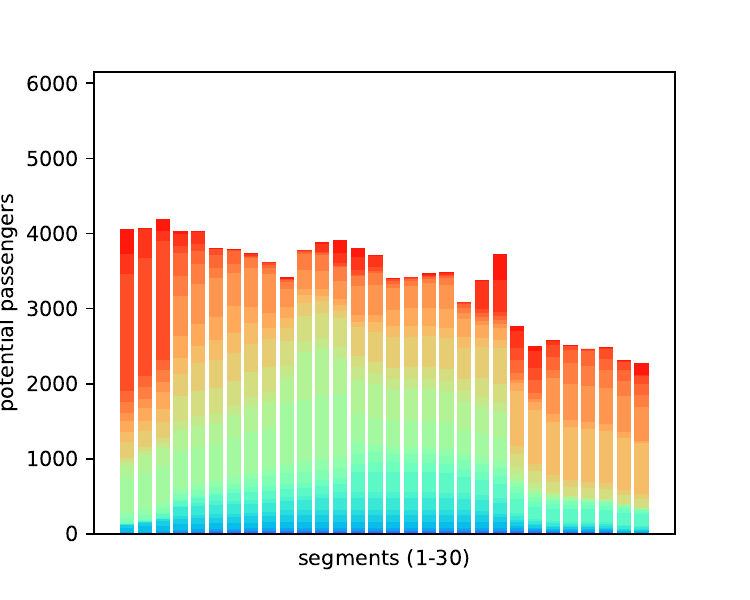}
		\caption{Alternative 3}
	\end{subfigure}
	\begin{subfigure}{0.32 \textwidth}
		\centering
		\includegraphics[width=\textwidth]{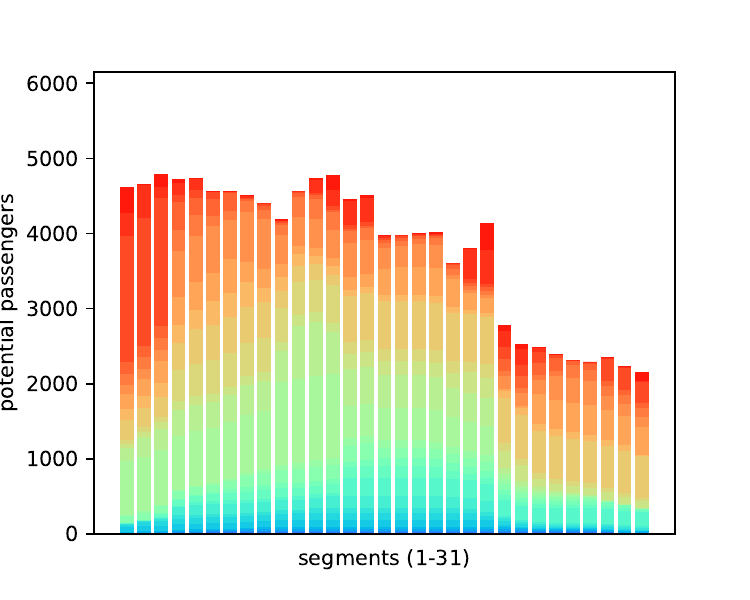}
		\caption{Alternative 4}
	\end{subfigure}
	\begin{subfigure}{0.32 \textwidth}
		\centering
		\includegraphics[width=\textwidth]{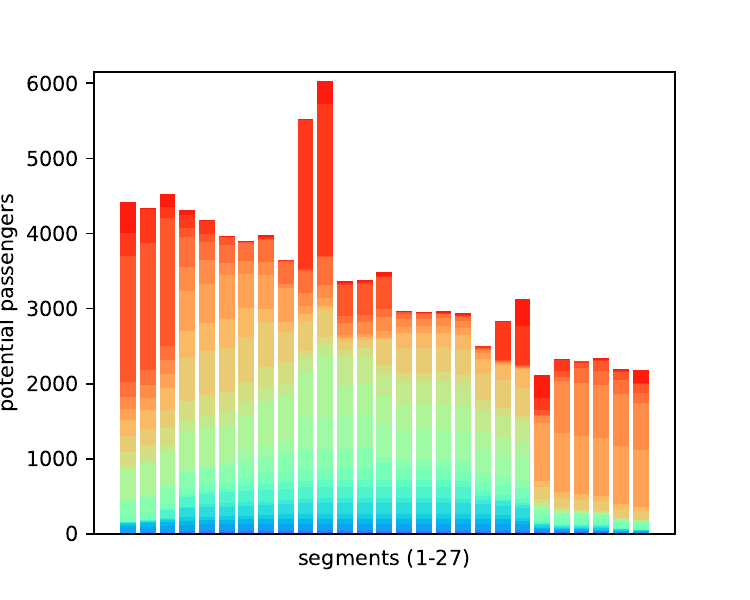}
		\caption{Alternative 5}
	\end{subfigure}
	\begin{subfigure}{0.32 \textwidth}
		\centering
		\includegraphics[width=0.85\textwidth]{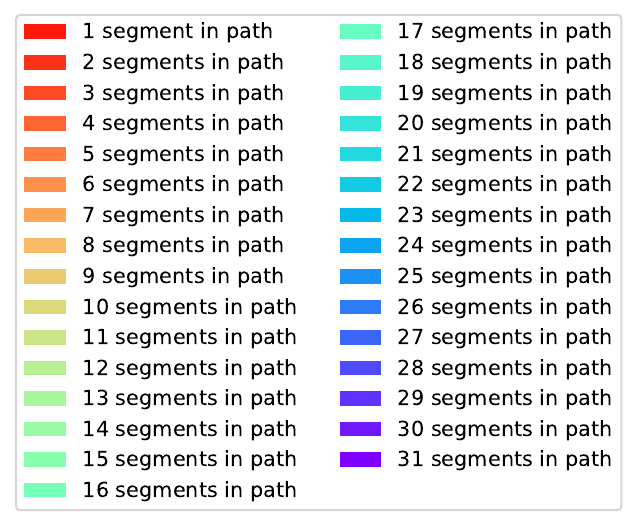}
		\caption{Legend}
	\end{subfigure}
	
	\caption{Load profiles for the five route alternatives. The horizontal axis contains each segment of the considered route alternative from north (Lyngby St.) to south (Ishøj St.). Each bar represents the number of potential passengers using the respective segment, and the coloring of a bar depicts the total travel distance of passengers using that segment.}
	\label{fig: load profiles Movia}
\end{figure}

In addition, we constructed an estimate of the future OD matrix by combining the estimated passenger demand per station from the pre-assessment study with the current OD matrix on the existing bus line. 
Specifically, the distribution of the forecasted inflow of passengers over all possible destination stations was scaled by the current fractional relation between this station and the other stations on the 400S line. 
A customized mapping was built for OD pairs that did not exist yet on the 400S. 
The resulting load profiles were determined based on conversations with Movia.
The obtained load profiles are shown in \Cref{fig: load profiles Movia}, where the height again indicates the load per segment, and the coloring indicates the length of the boarded passengers' paths.
We assume that a fixed percentage of each OD pair can be attributed to passengers newly attracted by infrastructure improvements.

We consider two potential budget splits between the seven municipalities that are willing to invest in the line.
These are the cost-based and the passenger-based budget splits, $\beta \in \{\textrm{COST}, \textrm{PASS}\}$, as described in \Cref{tab:param}. 
Because the actual costs and the number of passengers per municipality vary strongly, considering the \textrm{EQUAL} budget split is unrealistic. 
Also, the impact of the number of allowed BRT components is evaluated according to the parameters included in \Cref{tab:param}.
There are thus $3 \cdot 4 = 12$ instances per route alternative, giving a total of 60 instances for the case study.
These instances are available at \url{https://doi.org/10.11583/DTU.23664069}.

\subsection{Analysis of Pareto Fronts}
\label{sec:Copenhagen-case-study}
\begin{figure}[htbp]
	\centering
	\begin{subfigure}[b]{0.475\textwidth}
		\centering
		\includegraphics[width=0.95\textwidth]{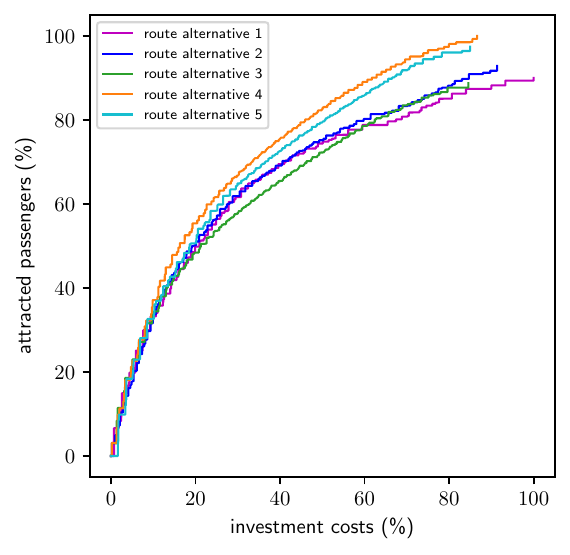}
		\caption{$|M| = 1$, \Linear{}}
		\label{fig:pareto_lin_SOC}
	\end{subfigure}
	\begin{subfigure}[b]{0.475\textwidth}
		\centering
		\includegraphics[width=0.95\textwidth]{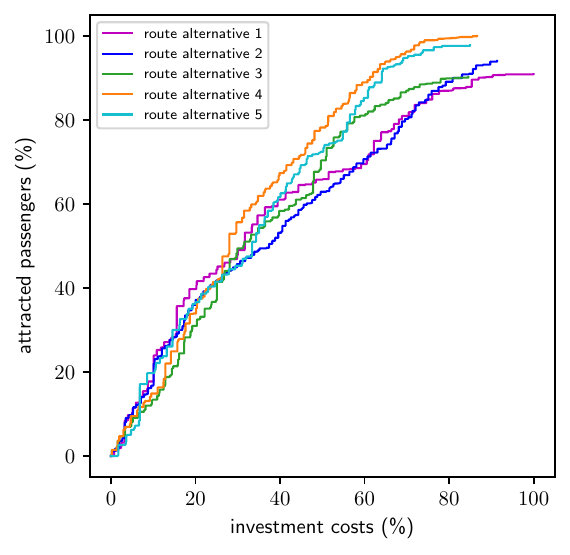}
		\caption{$|M| = 1$, \MinImprov{}}
		\label{fig:pareto_min_SOC}
	\end{subfigure}
	\begin{subfigure}[b]{0.475\textwidth}
		\centering
		\includegraphics[width=0.95\textwidth]{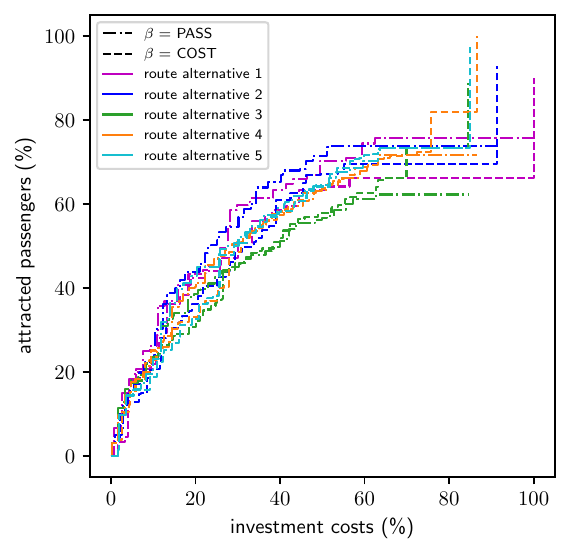}
		\caption{$|M| = 7$, \Linear{}}
		\label{fig:pareto_lin_BRT}
	\end{subfigure}
	\begin{subfigure}[b]{0.475\textwidth}
		\centering
		\includegraphics[width=0.95\textwidth]{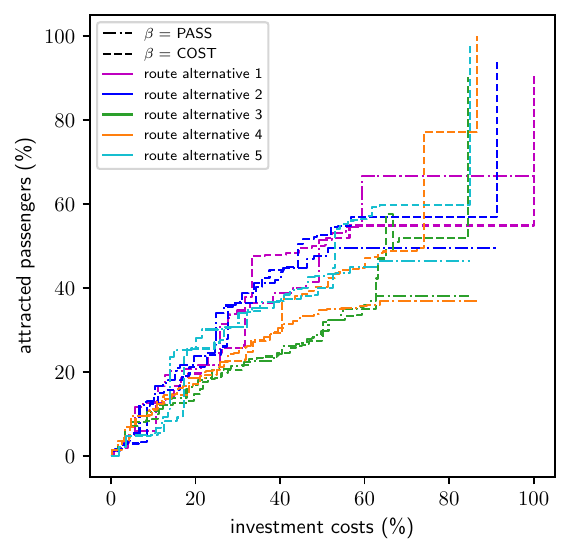}
		\caption{$|M| = 7$, \MinImprov{}}
		\label{fig:pareto_min_BRT}
	\end{subfigure}
	
	\caption{Comparing investment costs and attracted passengers for the five different route alternatives for $Z=\infty$. The investment costs are given as a percentage of the costs for upgrading all segments of the most expensive route alternative, and the attracted passengers are given as a percentage of the maximum number of potential passengers over all route alternatives.}
	\label{fig:case_study_pareto}
\end{figure}

We now look at the results of our experiments, where our aim is to analyze and compare the investment trade-offs for the five BRT route alternatives, taking into account the effect of the different passenger responses and budget splits over the municipalities.
Here, we use a similar computational set-up as described for the artificial instances in \Cref{sec: instances}.
Moreover, the investment budget is limited to the investment costs for upgrading all segments of the most expensive route alternative.
As a consequence, dependent on the budget split, not all municipalities may have sufficient budget to upgrade all of their segments.
The resulting Pareto fronts for the two passenger responses, with and without municipalities, evaluated regarding the investment costs are given in \Cref{fig:case_study_pareto} for the setting without a \ConConstr{} ($Z = \infty$).
Here, graphs in the top row provide the results when there is a global decision maker ($|M| = 1$), and graphs in the bottom row are for the case with municipalities ($|M| = 7$). 
Each graph indicates the investment costs as a percentage of the costs for upgrading all segments of the most expensive route alternative on the horizontal axis. 
The vertical axis indicates the attracted passengers relative to the maximum number of potential passengers over all route alternatives. 
This scaling on both axes allows to directly compare the route alternatives to each other.

The obtained Pareto plots show that many of the observations from the artificial results carry over to the Greater Copenhagen case study.
It can be seen that the number of attracted passengers is generally again higher for the \Linear{} passenger response than for the \MinImprov{} passenger response, except for investment levels that are above 75\% to 80\%, and that this effect is more pronounced when including the different municipalities.
Moreover, the introduction of municipality budgets has again a significant impact on the number of attracted passengers, especially under the passenger response \MinImprov{}.
However, especially apparent in these case study results is the ability of the budget split $\beta = \textrm{COST}$ to achieve a significantly higher number of passengers at higher investment levels.
This effect can be attributed to the presence of segments with very high upgrade costs, which are hard to upgrade for municipalities when they are not awarded a share that is in line with these upgrade costs.

When focusing on the comparison of the route alternatives, \Cref{fig:case_study_pareto} shows that there is not a universal ordering of the route alternatives respective to the number of attracted passengers.
Instead, this ordering depends on both the investment level and the passenger response.
For example, it can be seen that route alternatives 4 and 5 lead to the largest number of attracted passengers for middle to high investment levels under both passenger responses for $|M| = 1$, which can be explained by the higher total passenger potential for these alternatives.
However, for $\vert M \vert =7$ and at an investment level between 30\% and 70\%, the route alternatives 1 and 2 yield the highest numbers of attracted passengers for both passenger responses.
For low investment levels, the numbers of attracted passengers deviate less between the route alternatives, but it depends on the precise investment level, which route alternative is best.

Our results thus show the importance of obtaining knowledge about the passenger response and the willingness of municipalities to invest before a final route alternative is chosen for the BRT line.

\begin{figure}[htbp]
	\centering
	\begin{subfigure}[b]{\textwidth}
		\centering
		\includegraphics[width=0.95\textwidth]{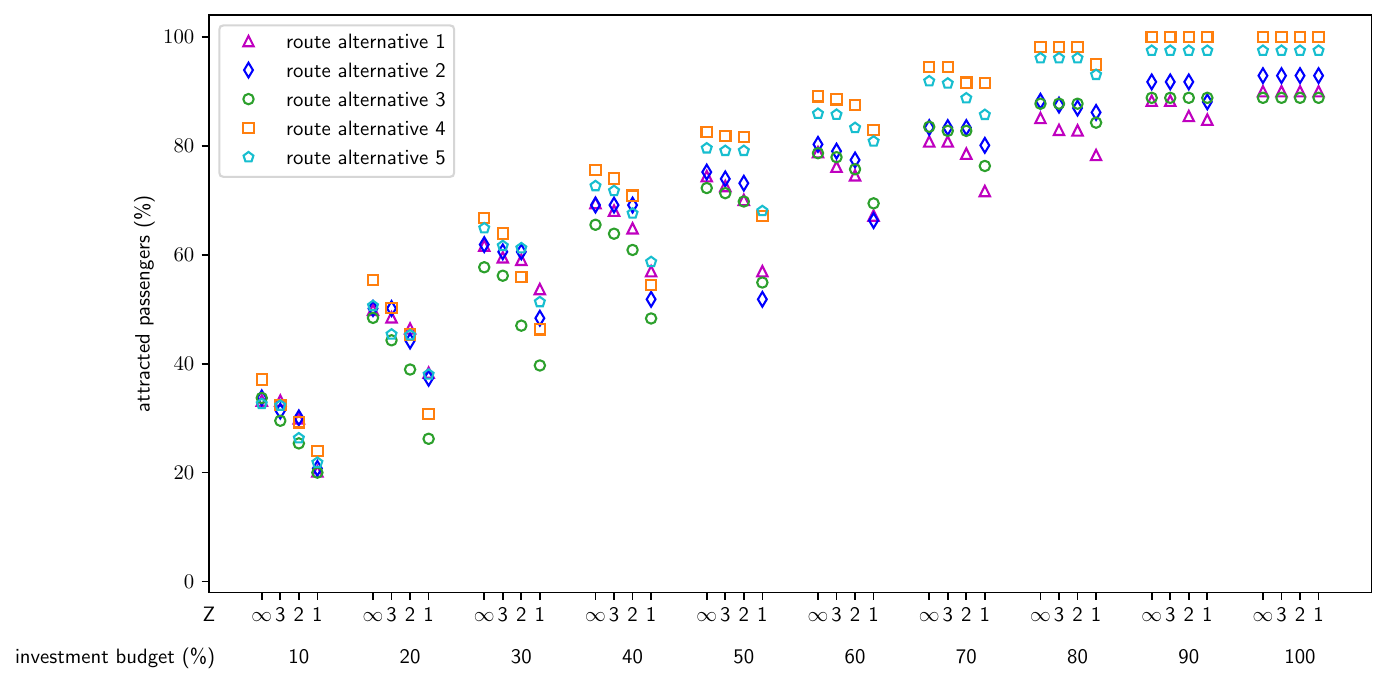}
		\caption{\Linear{}}
		\label{fig:mov_max_component_lin}
	\end{subfigure}
	\begin{subfigure}[b]{\textwidth}
		\centering
		\includegraphics[width=0.95\textwidth]{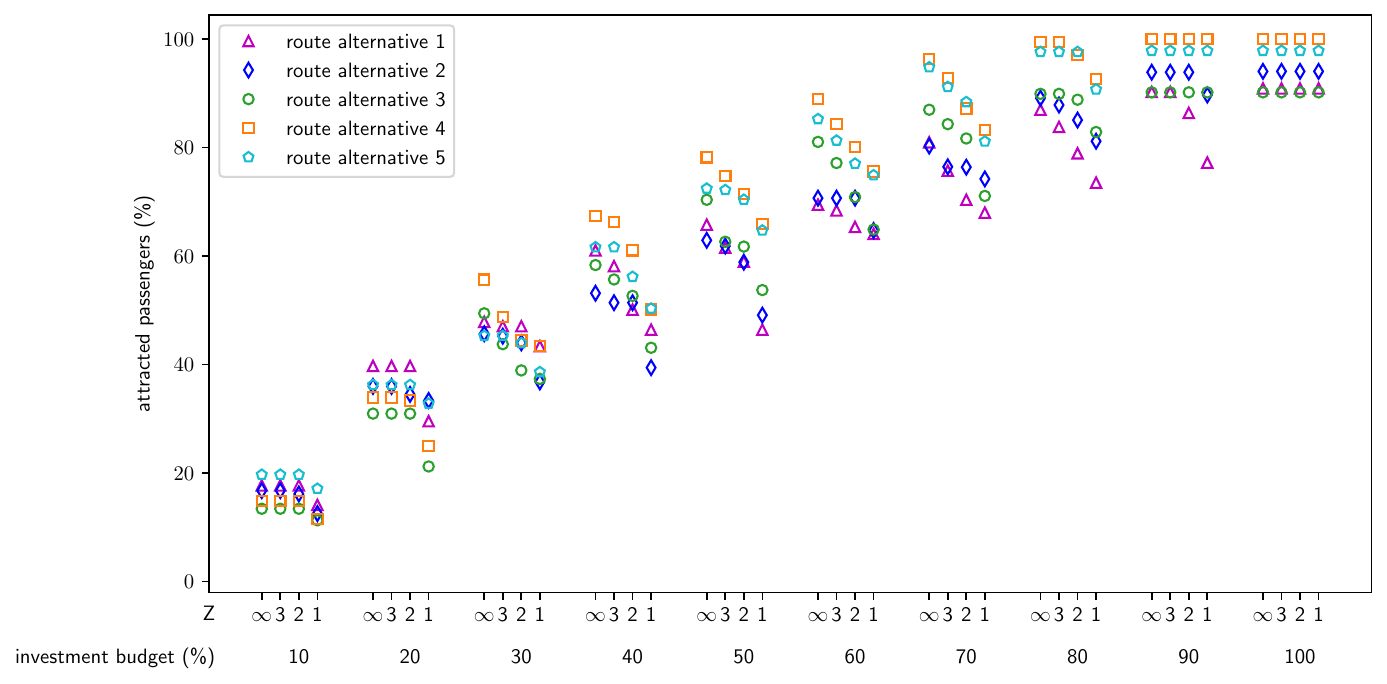}
		\caption{\MinImprov{}}
		\label{fig:mov_max_component_min}
	\end{subfigure}
	
	\caption{Influence of the upper bound on the number of BRT components on the percentage of attracted passengers for the five route alternatives. This figure condenses Pareto plots to ten budget ranges and shows the solution with the highest number of attracted passengers within each range. The investment budget is given as a percentage of the costs for upgrading all segments of the most expensive route alternative, and the attracted passengers are given as a percentage of the maximum number of potential passengers over all route alternatives.}
	\label{fig:mov_max_component}
\end{figure}

\paragraph{Influence of the Number of Components}
It remains to analyze the impact of the \ConConstr{} for the Greater Copenhagen case study.
This effect is depicted in \Cref{fig:mov_max_component}, which analyzes the effect of the number of allowed BRT components $Z$ on the number of attracted passengers for each of the five route alternatives and for both passenger responses.
These results are computed for the setting of a global decision maker ($|M| = 1$).
Moreover, to make the impact of the \ConConstr{} more visible, this figure condenses the Pareto plots to ten budget ranges and shows the solution with the highest number of attracted passengers within each range.

As for the artificial instances, \Cref{fig:mov_max_component} shows that restricting the number of components leads to a reduced number of attracted passengers for all route alternatives.
This effect is strongest for investment levels that are closer to the middle and lower end. 
By design, no effect can be seen for the highest investment level, as all segments will be upgraded.
Comparing the \Linear{} and \MinImprov{} passenger responses, an interesting difference is that the impact of restricting the number of components is stronger for the very low investment levels for \Linear{}.
In addition, it can be seen that it is again especially the restriction to a single component that leads to a strong reduction in passengers. 
For these instances, there is a difference between allowing 2, 3 or arbitrarily many BRT components for most available budget levels, although the solution for at most 3 BRT components comes close to that of allowing arbitrarily many BRT components.

\section{Conclusion}
\label{sec:conclusion}
We studied the bi-objective \BRTproblem{}, which focuses on determining the set of segments to be upgraded for a BRT line such to balance the number of attracted passengers and the investment budget.
Municipalities are considered in this problem through separate municipality budgets.
Moreover, this problem allows the restriction of the number of upgraded connected components to prevent frequent switching between upgraded and non-upgraded segments.
Additionally, we considered two passenger responses to upgrades: a linear and a threshold relation.

We developed a bi-objective mixed-integer linear programming formulation for the \BRTproblem{} and an algorithm based on the $\epsilon$-constraint method to find the complete set of non-dominated points.
We proved that the number of non-dominated points grows exponentially in general but identified special cases in which the problem becomes tractable.
Similarly, we showed that the subproblems that are solved within the $\epsilon$-constraint-based algorithm are generally NP-hard but allow polynomially solvable special cases.

Our numerical experiments for artificial instances and the Greater Copenhagen case study analyzed the impact of the passenger response, the separate municipality budgets and the \ConConstr{}.
The main findings indicate that splitting the budget over municipalities directly reduces the number of attracted passengers, as does the requirement to have only one BRT component.
However, as soon as two or three BRT components are allowed, the impact is far smaller.
Regarding the artificial instances, for investment costs below 75\% of the total costs for upgrading all segments, the linear passenger response indicates higher numbers of attracted passengers than the threshold passenger response.
For higher investment costs, it turns, but the values are quite close.
The Greater Copenhagen case study confirmed many of these observations, showing that they translate to real-life instances.
Further, the Greater Copenhagen case study showed that the ranking of the route alternatives is highly dependent on both the passenger response and the available investment budget.
Hence, obtaining a good estimate on how passengers respond to the upgrades and on the extent to which municipalities are willing to invest is crucial for selecting the best route alternative.

In this paper, we have analyzed the two extreme cases of a linear (\Linear{}) and threshold (\MinImprov{}) passenger response to upgrades.
Considering mixes of these two passenger response functions would, therefore, be a natural next step.
Such a mix could, e.g., be a piecewise linear response function to upgrades, where the impact of an upgrade depends on the overall extent to which upgrades are realized.
This would include the special case where the number of passengers grows linearly as soon as a certain threshold of infrastructure improvements is reached.
Note that piecewise linear objective functions can be easily integrated into the algorithmic approach suggested in this paper, meaning that the suggested $\epsilon$-constraint-based algorithm can still be used to find the complete Pareto front.

In addition, the context of this paper assumes an interest in a social optimum, provided a certain investment level per municipality.
As discussed, such a setting could occur when there is a third party, such as a transport agency, that makes suggestions to the municipalities as to which segments should be upgraded.
One could imagine that in certain cases, the interest for a social optimum is challenged by the individual interests of municipalities, for example, when not all municipalities benefit equally from the investments to be made. 
It would therefore be interesting for future work to focus on a game-theoretic setting that models this competition, possibly within the context of a central government searching for the best subsidy scheme to attract the most passengers at minimum budget subject to the internal competition between municipalities.

Another interesting direction of future work could be considering the investment problem in a network context instead of for a single line, either by including the determination of the route of the BRT line or by considering that other lines could (be rerouted to) profit from the upgraded BRT infrastructure as well.
In a network setting, other models and aspects for passenger behavior could be considered, e.g., including the travel time \citep{schiewe2020periodic} and fares \citep{Schoebel2022} as well as route and mode choice.
Also, the inclusion of operating constraints considering load profiles, e.g., in the setting of self-driving minibusses with innovative operating modes \citep{gkiotsalitis2022frequency} could be an interesting direction.

\section*{Acknowledgement}
This work was supported by the European Union’s Horizon 2020 research and innovation programme [Grant 875022], the Federal Ministry of Education and Research [Project 01UV2152B], and Innovationsfonden [0205-00002B] under the project sEAmless SustaInable EveRyday urban mobility (EASIER). Moreover, we would like to thank Region H [0205-00005B] and Movia for providing insight into the planning process of the BRT system and for the provision of data.

\bibliographystyle{abbrvnat}
\bibliography{easier_papers}

\clearpage
\section*{Appendix: Additional Instance Details and Evaluations}
\label{sec:appendix}

\begin{figure}[h!t]
	\centering
	\begin{subfigure}[b]{0.85\textwidth}
		\centering
		\includegraphics[width=0.95\textwidth]{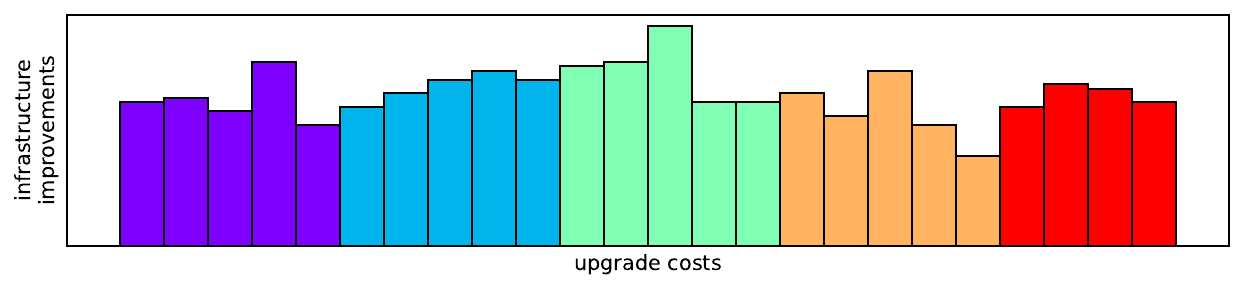}
		\caption{UNIT}
	\end{subfigure}
	\begin{subfigure}[b]{0.85\textwidth}
		\centering
		\includegraphics[width=0.95\textwidth]{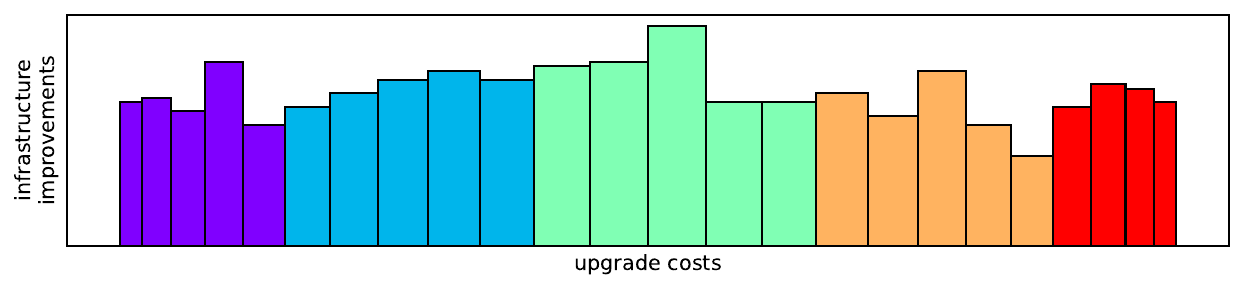}
		\caption{MIDDLE}
	\end{subfigure}
	\begin{subfigure}[b]{0.85\textwidth}
		\centering
		\includegraphics[width=0.95\textwidth]{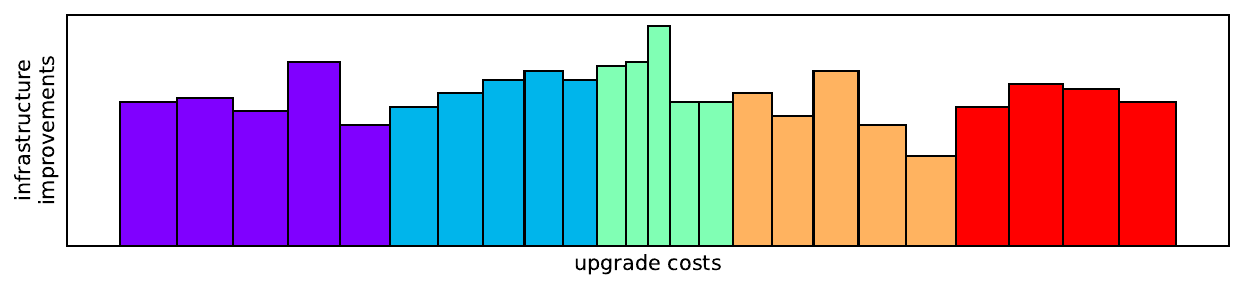}
		\caption{ENDS}
	\end{subfigure}
	\caption{Cost patterns and infrastructure improvements per segment for the artificial instances. Each bar represents a segment. The width of a bar represents the upgrade costs while the height reflects the infrastructure improvements. The colors indicate to which municipality a segment belongs.}
	\label{fig:costs}
\end{figure}

\begin{figure}[h!t]
	\begin{subfigure}[b]{0.99\textwidth}
		\centering
		\resizebox{\textwidth}{!}{
			\begin{tikzpicture}[thick, main/.style = {draw, circle, fill=black, inner sep=1pt,minimum size=.4cm}, small/.style = {draw, circle, inner sep=1pt,minimum size=.2cm}]
				\foreach \i in {1,...,24}
				\node[small] (x{\i}) at (\i,0){};
				\foreach \i in {1,...,23}
				\draw (x{\i}) -- (x{\the\numexpr \i+1});
			\end{tikzpicture}
		}
		\caption{\EvenDemand}
	\end{subfigure}
	
	\vspace{1em}
	\begin{subfigure}[b]{0.99\textwidth}
		\centering
		\resizebox{\textwidth}{!}{
			\begin{tikzpicture}[thick, main/.style = {draw, circle, fill=black, inner sep=1pt,minimum size=.4cm}, small/.style = {draw, circle, inner sep=1pt,minimum size=.2cm}]
				\foreach \i in {1,...,25}
				\ifthenelse{\i =1 \OR \i=9 \OR \i=17}{\node[main] (x{\i}) at (\i,0) {}}{\node[small] (x{\i}) at (\i,0){}};
				\foreach \i in {1,...,24}
				\draw (x{\i}) -- (x{\the\numexpr \i+1});
			\end{tikzpicture}
		}
		\caption{\LargeStationsDemand}
	\end{subfigure}
	
	\vspace{1em}
	\begin{subfigure}[b]{0.99\textwidth}
		\centering
		\resizebox{\textwidth}{!}{
			\begin{tikzpicture}[thick, main/.style = {draw, circle, fill=black, inner sep=1pt,minimum size=.4cm}, small/.style = {draw, circle, inner sep=1pt,minimum size=.2cm}]
				\foreach \i in {1,...,25}
				\ifthenelse{\i =1 \OR \i=25}{\node[main] (x{\i}) at (\i,0) {}}{\node[small] (x{\i}) at (\i,0){}};
				\foreach \i in {1,...,24}
				\draw (x{\i}) -- (x{\the\numexpr \i+1});
			\end{tikzpicture}
		}
		\caption{\EndStationsDemand}
	\end{subfigure}
	\caption{Line graph of the artificial instances. Stations with high demand in the demand patterns \LargeStationsDemand{} and \EndStationsDemand{} are marked with filled black nodes. There are no large stations in demand pattern \EvenDemand{}.}
	\label{fig:graphs}
\end{figure}
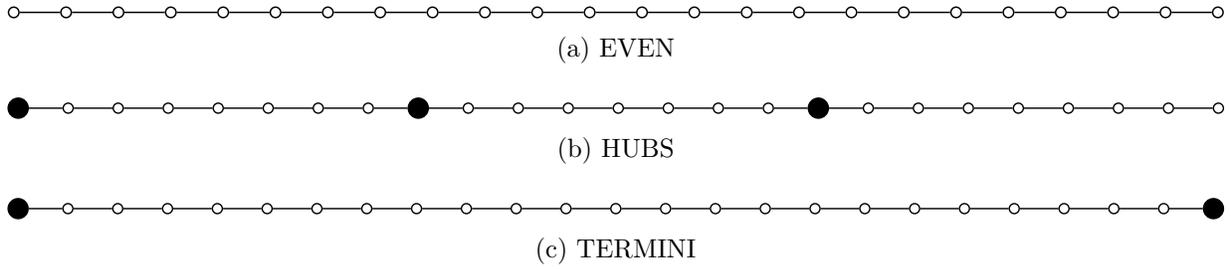

\begin{figure}[h!t]
	\begin{subfigure}[b]{0.5\textwidth}
		\centering
		\includegraphics[width=\textwidth]{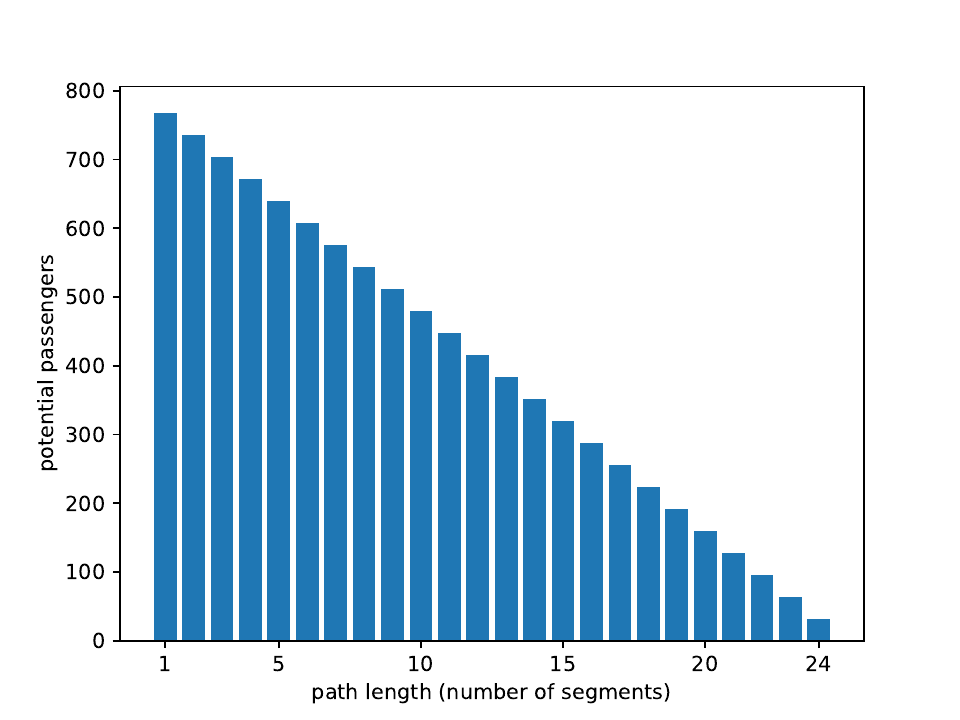}
		\caption{\EvenDemand}
	\end{subfigure}
	\begin{subfigure}[b]{0.5\textwidth}
		\centering
		\includegraphics[width=\textwidth]{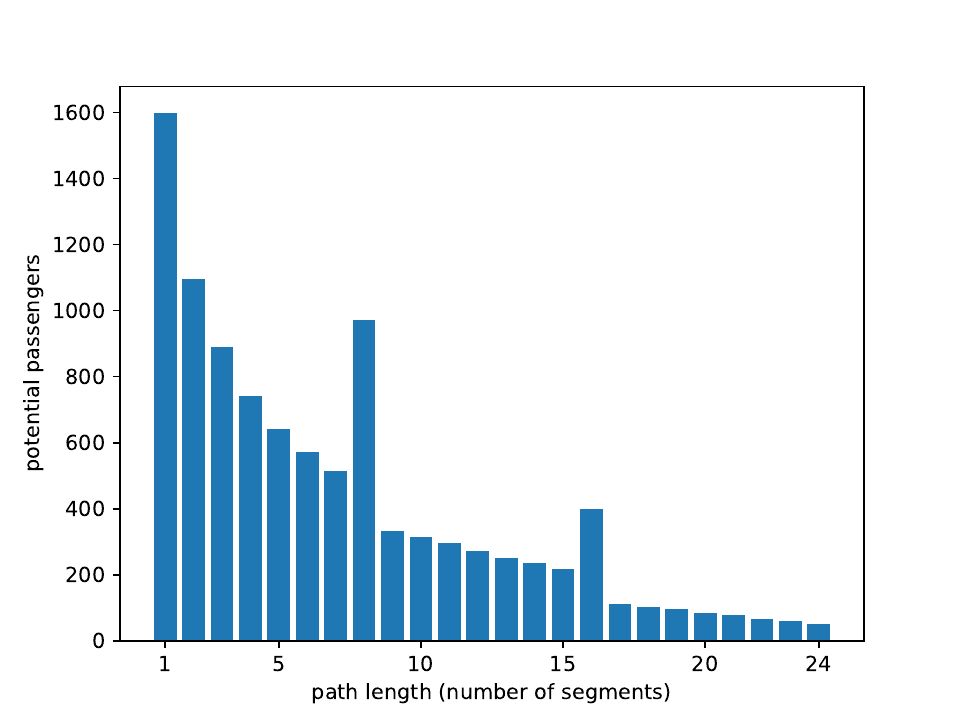}
		\caption{\LargeStationsDemand}
	\end{subfigure}
	\begin{subfigure}[b]{0.5\textwidth}
		\centering
		\includegraphics[width=\textwidth]{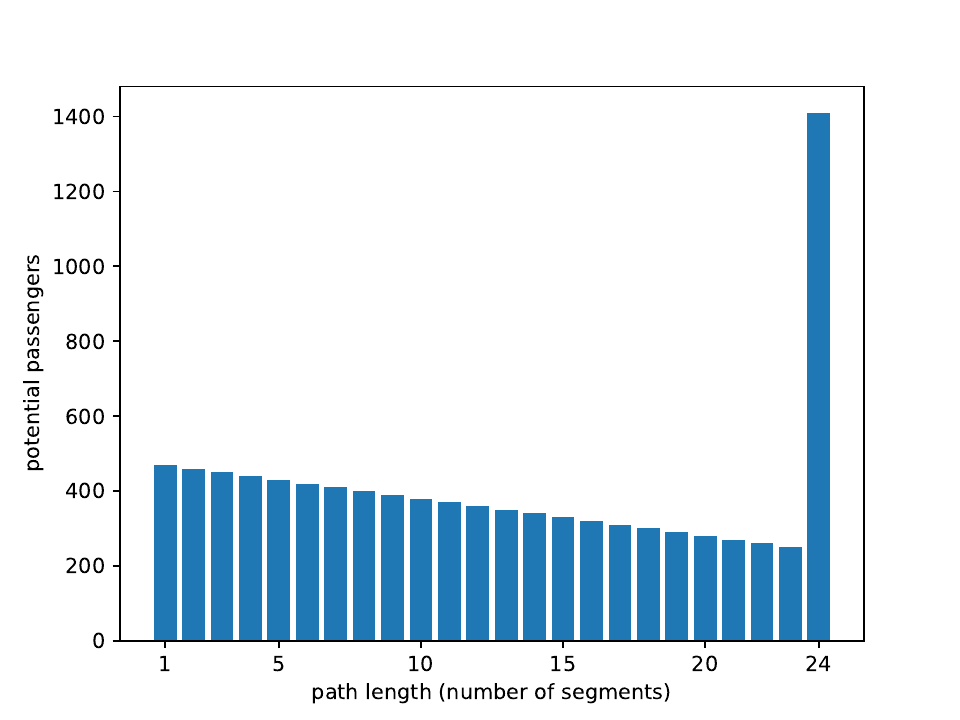}
		\caption{\EndStationsDemand}
	\end{subfigure}
	\caption{Histogram of the travel distances of passengers. The height of a bar gives the demand of passengers traveling for a certain number of segments.}
	\label{fig:hist-distances}
\end{figure}

\begin{figure}[ht]
	\centering
	\begin{subfigure}[b]{0.475\textwidth}
		\centering
		\includegraphics[width=0.95\textwidth]{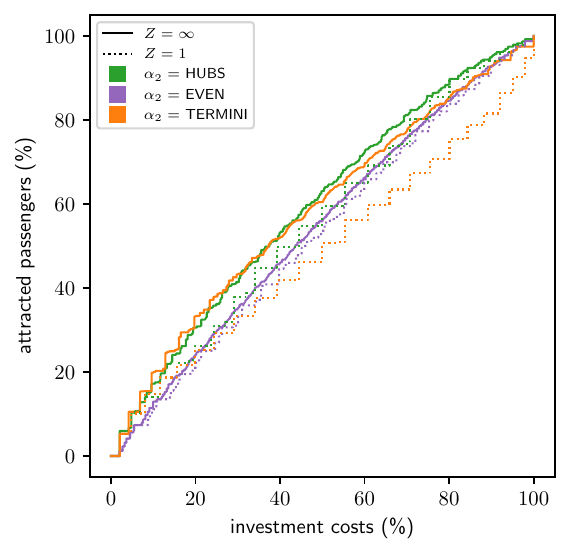}
		\caption{\Linear{}, $\alpha_1 = \textrm{MIDDLE}$.}
	\end{subfigure}
	\begin{subfigure}[b]{0.475\textwidth}
		\centering
		\includegraphics[width=0.95\textwidth]{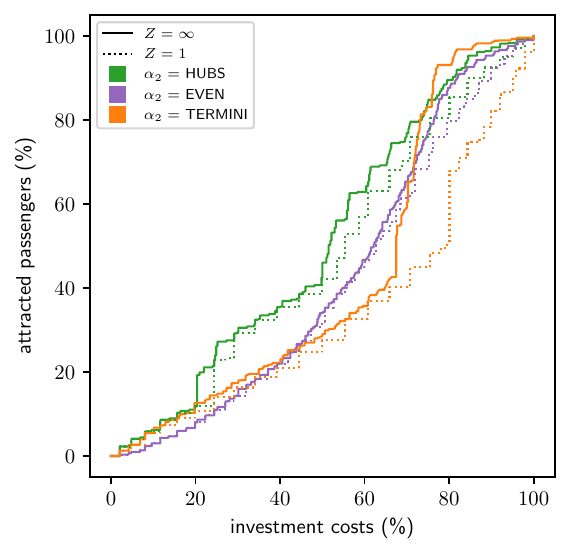}
		\caption{\MinImprov{}, $\alpha_1 = \textrm{MIDDLE}$.}
	\end{subfigure}
	\begin{subfigure}[b]{0.475\textwidth}
		\centering
		\includegraphics[width=0.95\textwidth]{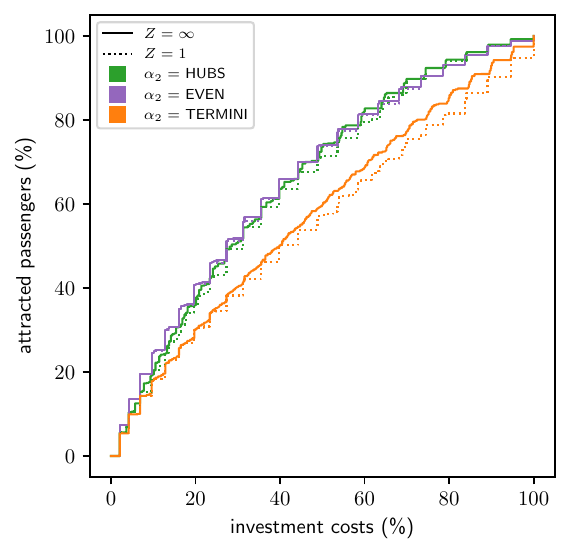}
		\caption{\Linear{}, $\alpha_1 = \textrm{ENDS}$.}
	\end{subfigure}
	\begin{subfigure}[b]{0.475\textwidth}
		\centering
		\includegraphics[width=0.95\textwidth]{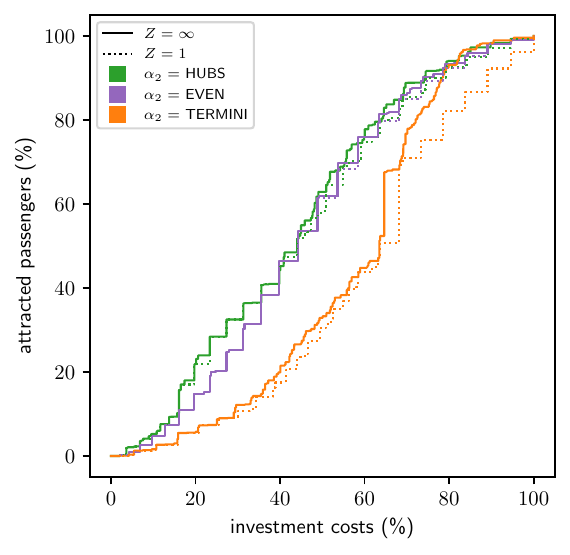}
		\caption{\MinImprov{}, $\alpha_1 = \textrm{ENDS}$.}
	\end{subfigure}
	\caption{Evaluation of the non-dominated points of $\SOCtwo{\star}{Z\geq 1}$ for artificial instances with cost pattern $\alpha_1 \in \{ \textrm{ENDS},\textrm{MIDDLE} \}$ and $Z \in \{1,\infty\}$. Both attracted passengers and investment costs are given as percentage of the total number of potential passengers and costs for upgrading all segments, respectively.}
	\label{fig:pareto_25_cost_pattern}
\end{figure}

\begin{figure}[h!t]
	\centering
	\begin{subfigure}[b]{\textwidth}
		\centering
		\includegraphics[width=0.83\textwidth]{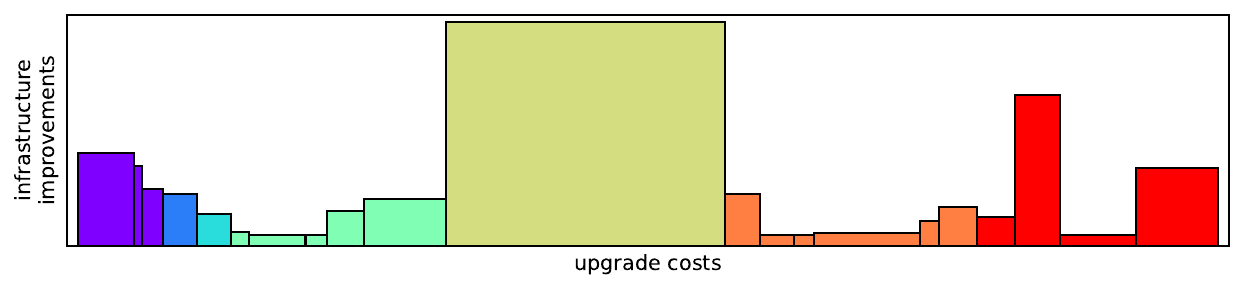}
		\caption{Alternative 1}
	\end{subfigure}
	\begin{subfigure}[b]{\textwidth}
		\centering
		\includegraphics[width=0.83\textwidth]{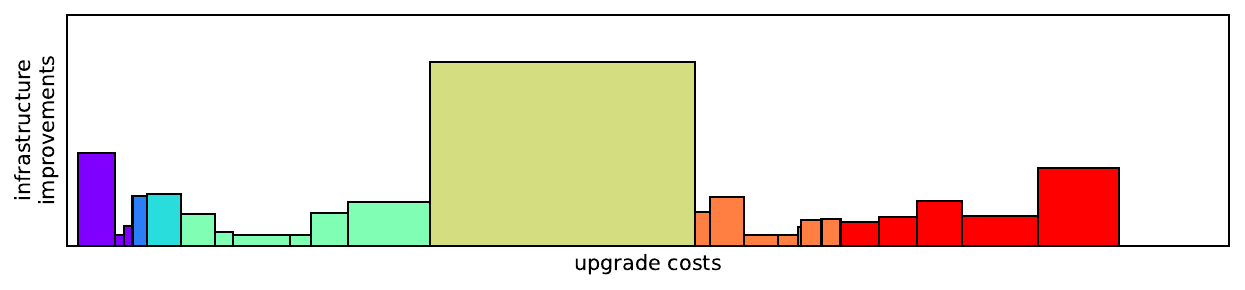}
		\caption{Alternative 2}
	\end{subfigure}
	\begin{subfigure}[b]{\textwidth}
		\centering
		\includegraphics[width=0.83\textwidth]{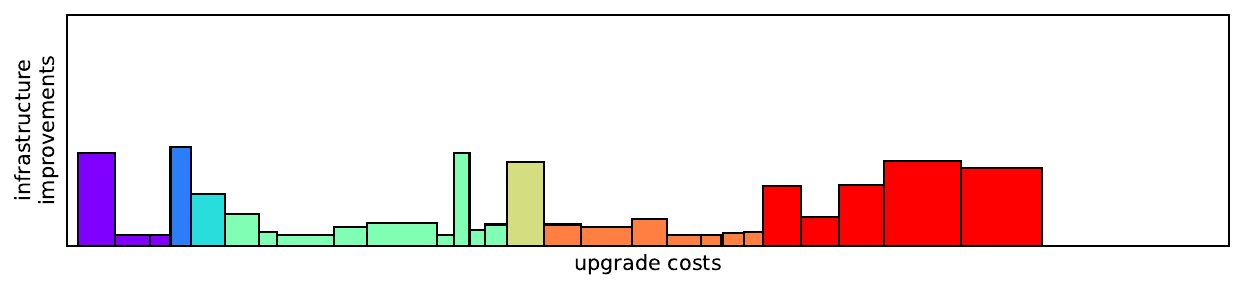}
		\caption{Alternative 3}
	\end{subfigure}
	\begin{subfigure}[b]{\textwidth}
		\centering
		\includegraphics[width=0.83\textwidth]{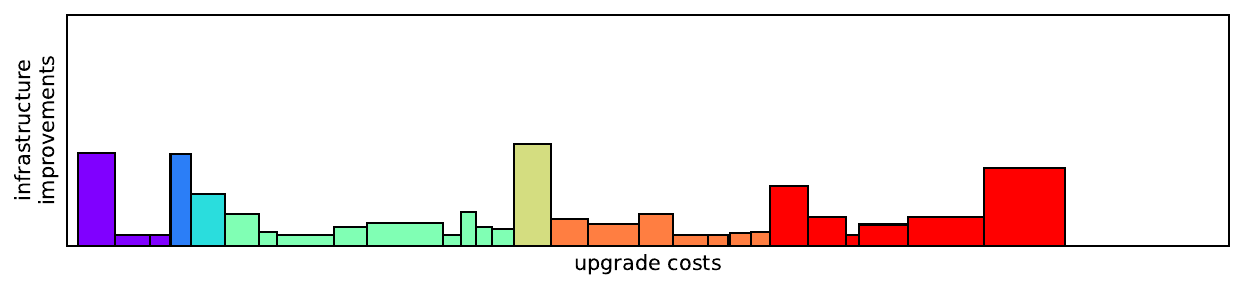}
		\caption{Alternative 4}
	\end{subfigure}
	\begin{subfigure}[b]{\textwidth}
		\centering
		\includegraphics[width=0.83\textwidth]{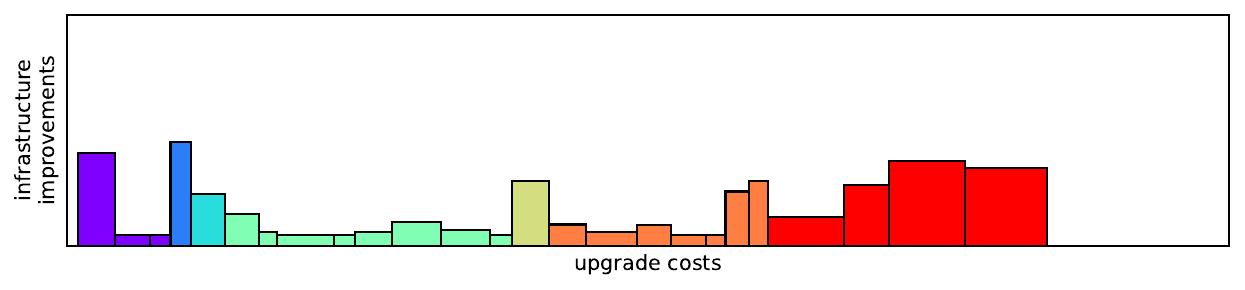}
		\caption{Alternative 5}
	\end{subfigure}
	\caption{Cost patterns and infrastructure improvements per segment for the five route alternatives from north (Aldershvilevej) to south (Ishøj St.). Each bar represents a segment. The width of a bar represents the upgrade costs while the height reflects the infrastructure improvements. The colors indicate to which municipality a segment belongs. Note that the two non-upgradable segments in Lyngby municipality are excluded.}
	\label{fig:costs-movia}
\end{figure}

In the appendix, additional information about the numerical experiments and the Greater Copenhagen case study is provided.
\Cref{fig:costs} and \Cref{fig:costs-movia} show the upgrade costs and the infrastructure improvements per segment as well as which segments belong to the same municipality for the artificial instances and the case study, respectively.
\Cref{fig:graphs} depicts the graphs of the artificial instances, marking the locations of stations with high demand within the demand patterns \LargeStationsDemand{} and \EndStationsDemand{}.
\Cref{fig:hist-distances} shows the corresponding histograms of the travel distances of the passengers.
Moreover, as a supplement to \Cref{fig:pareto_25_unit}, \Cref{fig:pareto_25_cost_pattern} shows the evaluation of the non-dominated points of \SOCtwo{$\star$}{$Z\geq 1$} for the artificial instances with the cost patterns ENDS and MIDDLE.
Further plots that depict which segments are upgraded at certain investment budget levels are provided at \url{https://doi.org/10.11583/DTU.c.6805470}.

\end{document}